\input amstex
\input amsppt.sty
\magnification=\magstep1
\hsize=34truecc
 \vsize=22.5truecm
\baselineskip=14truept
\NoBlackBoxes
\def\q{\quad}
\def\qq{\qquad}
\def\mod#1{\ (\text{\rm mod}\ #1)}
\def\t{\text}
\def\qtq#1{\q\t{#1}\q}
\def\f{\frac}
\def\e{\equiv}
\def\a{\alpha}

\def\ep{\varepsilon}

\def\sls#1#2{(\f{#1}{#2})}
 \def\ls#1#2{\big(\f{#1}{#2}\big)}
\def\Ls#1#2{\Big(\f{#1}{#2}\Big)}
\def\sqs#1#2{(\f{#1}{#2})_4}
\def\qs#1#2{\big(\f{#1}{#2}\big)_4}
\def\Qs#1#2{\Big(\f{#1}{#2}\Big)_4}

\let \pro=\proclaim
\let \endpro=\endproclaim
\topmatter
\title\nofrills Congruences for $q^{[p/8]}\pmod p$ II
\endtitle
\author ZHI-Hong Sun\endauthor
\affil School of Mathematical Sciences, Huaiyin Normal University,
\\ Huaian, Jiangsu 223001, PR China
\\ E-mail: zhihongsun$\@$yahoo.com
\\ Homepage: http://www.hytc.edu.cn/xsjl/szh
\endaffil
 \nologo \NoRunningHeads
 \abstract{Let $\Bbb Z$ be the set of integers, and let $p$
be a prime of the form $4k+1$. Suppose $q\in\Bbb Z$, $2\nmid q$,
$p\nmid q$, $p=c^2+d^2$, $c,d\in\Bbb Z$ and $c\equiv 1\pmod 4$. In
this paper we continue to discuss congruences for $q^{[p/8]}\pmod p$
and present new reciprocity laws, but we assume $4p=x^2+qy^2$ or
$p=x^2+2qy^2$, where $[\cdot]$ is the greatest integer function and
$x,y\in\Bbb Z$.

\par\q
}
 \endabstract
 \footnote"" {The author is
supported by the National Natural Sciences Foundation of China (No.
11371163).}
\endtopmatter
\document
\subheading{1. Introduction}
\par Let $\Bbb Z$ be the set of integers, $i=\sqrt{-1}$ and $\Bbb Z[i]=\{a+bi\mid a,b\in\Bbb
Z\}$. For any positive odd number $m$ and $a\in\Bbb Z$ let $(\f am)$
be the (quadratic) Jacobi symbol. (We also assume $(\f a1)=1$.) For
our convenience we also define $(\f a{-m})=(\f am)$. Then for any
two odd numbers $m$ and $n$ with $m>0$ or $n>0$ we have the
following general quadratic reciprocity law: $\sls
mn=(-1)^{\f{m-1}2\cdot \f{n-1}2}\sls nm$.
\par  For $a,b,c,d\in\Bbb Z$ with $2\nmid c$ and $2\mid d$, one can define the
quartic Jacobi symbol $\qs{a+bi}{c+di}$ as in [S1,S2,S4]. From [IR]
we know that
$\overline{\qs{a+bi}{c+di}}=\qs{a-bi}{c-di}=\qs{a+bi}{c+di}^{-1}$,
where $\bar x$ means the complex conjugate of $x$.
 In Section 2 we list main properties of the quartic Jacobi
symbol. See also [IR], [BEW] and [S4]. For the history of quartic
reciprocity laws, see [Lem].

\par Let $p$ be a prime of the form $4k+1$, $q\in\Bbb Z$, $2\nmid q$ and
$p\nmid q$. Suppose that $p=c^2+d^2=x^2+qy^2$, $c,d,x,y\in\Bbb Z$,
$c\e 1\mod 4$, $d=2^rd_0$ and $d_0\e 1\mod 4$. Assume that
$(c,x+d)=1$ or $(d_0,x+c)=1$, where $(m,n)$ is the greatest common
divisor of $m$ and $n$. In [S5], using the quartic reciprocity law
the author deduced some congruences for $q^{[p/8]}\mod p$ in terms
of $c,d,x$ and $y$, where $[\cdot]$ is the greatest integer
function. For example, if $p$ and $q$ are primes of the form $8k+1$,
$p=c^2+d^2=x^2+qy^2$, $q=a^2+b^2$, $a,b,c,d,x,y\in\Bbb Z$, $c\e
1\mod 4$, $d=2^rd_0$ and $d_0\e 1\mod 4$, then
$$q^{\f{p-1}8}\e (-1)^{\f d4+\f{xy}4}\Ls dc^m\mod p
\iff \Ls{ac+bd}{ac-bd}^{\f{q-1}8}\e \Ls ba^m\mod q .$$ If $p\e
1,9,25\mod{56}$ is a prime and so $p=c^2+d^2=x^2+7y^2$ for some
$c,d,x,y\in\Bbb Z$, by [S5, Corollary 4.2] we have
$$7^{\f{p-1}8}\e \cases-(-1)^{\f y4}\mod p&\t{if $7\mid c$,}
\\(-1)^{\f y4}\mod p&\t{if $7\mid d$,}
\\\pm(-1)^{\f y4}\f cd\mod p&\t{if $c\e \pm d\mod 7$.}
\endcases$$
\par Let $p$
be a prime of the form $4k+1$, $q\in\Bbb Z$, $2\nmid q$ and $p\nmid
q$. Suppose that $p=c^2+d^2$ with $c,d,x,y\in\Bbb Z$ and $c\e 1\mod
4$. In this paper we continue to discuss congruences for
$q^{[p/8]}\mod p$, but we assume $4p=x^2+qy^2$ or $p=x^2+2qy^2$.
Here are some typical results: \par $\star$ Let $p$ be a prime of
the form $8k+1$ and $\sls{-11}p=1$, $p=c^2+d^2$, $4p=x^2+11y^2$,
$c,d,x,y\in\Bbb Z$, $c\e 1\mod 4$ and $d=2^rd_0(2\nmid d_0)$. Assume
that $(c,x+2d)=1$ or $(d_0,x+2c)=1$. Then
$$(-11)^{\f {p-1}8}\e
\cases \pm(-1)^{[\f x4]}\mod p&\t{if $2\nmid x$ and $x\e \pm 4c,\pm
9c\mod{11}$,}
\\\pm (-1)^{[\f x4]}\f dc\mod p&\t{if $2\nmid x$
 and $x\e \pm 4d,\pm 9d\mod{11}$,}
\\\mp(-1)^{[\f x8]+\f y8}\mod p&\t{if $2\mid x$ and $x\e \pm 4c,\pm
9c\mod{11}$,}
\\\mp (-1)^{[\f x8]+\f y8}\f dc\mod p&\t{if $2\mid x$ and $x\e \pm 4d,\pm 9d\mod{11}$.}
\endcases$$
\par $\star$ Let $p$ and $q$ be primes such that $p\e 1\mod 8$, $q\e 7\mod 8$,
$p=c^2+d^2=x^2+2qy^2$, $c,d,x,y\in\Bbb Z$, $c\e 1\mod 4$, $d=2^rd_0$
and $d_0\e 1\mod 4$.  Assume $(c,x+d)=1$ or $(d_0,x+c)=1$. Then
 $$(-q)^{\f{p-1}8}\e \Ls dc^m\mod p\iff
\Ls{c-di}{c+di}^{\f{q+1}8}\e i^m\mod q.$$
\par $\star$ Let $p$
be a prime of the form $8k+1$ and $a\in\Bbb Z$ with $2\nmid a$.
Suppose that $p=c^2+d^2=x^2+(a^2+1)y^2$, $c,d,x,y\in\Bbb Z$, $c\e
1\mod 4$, $d=2^rd_0$, $d_0\e  1\mod 4$ and $4\mid y$. Assume
$(c,x+d)=1$ or $(d_0,x+c)=1$. Then
$$(a+\sqrt{a^2+1})^{\f{p-1}4}\e (-1)^{\f d4+\f y4}\mod p.$$
\par $\star$ Let $p\e 1,9\mod{40}$ be a
prime and so $p=c^2+d^2=x^2+10y^2$ with $c,d,x,y\in\Bbb Z$. Suppose
$c\e 1\mod 4$, $d=2^rd_0$ and $d_0\e 1\mod 4$. Assume $(c,x+d)=1$ or
$(d_0,x+c)=1$. Then
$$(-5)^{\f{p-1}8}\e\cases \pm (-1)^{\f d4+\f {x-1}2+\f y2}\mod
p&\t{if $\f{2c-d}x\e \pm 2\mod 5$,}
\\\pm (-1)^{\f d4+\f{x-1}2+\f y2}\f cd\mod p&\t{if $\f{2c-d}x\e \pm
1\mod 5$.}
\endcases$$
\par Throughout this paper, if $n\in\Bbb Z$, $2^{\a}\mid n$ and $2^{\a+1}\nmid
n$, we then write that $2^{\a}\ \Vert \ n$.

 \subheading{2. Basic lemmas}

 \pro{Lemma 2.1 ([S4, Proposition 2.1])} Let $a,b\in\Bbb Z$ with $2\nmid a$
and $2\mid b$. Then
$$\Qs i{a+bi}=i^{\f{a^2+b^2-1}4}
=(-1)^{\f{a^2-1}8}i^{(1-(-1)^{\f b2})/2}$$ and
$$\aligned\Qs{1+i}{a+bi}=\cases i^{((-1)^{\f{a-1}2}(a-b)-1)/4}&\t{if $4\mid b$,}
\\i^{\f{(-1)^{\f{a-1}2}(b-a)-1}4-1}&\t{if $2\ \Vert\ b$.}\endcases
\endaligned$$
\endpro
\pro{Lemma 2.2 ([S4, Proposition 2.2])} Let $a,b\in\Bbb Z$ with
$2\nmid a$ and $2\mid b$. Then
$$\Qs{-1}{a+bi}=(-1)^{\f b2}\qtq{and}\Qs
2{a+bi}=i^{(-1)^{\f{a-1}2}\f b2}.$$
\endpro
\pro{Lemma 2.3 ([S4, Proposition 2.3])} Let $a,b,c,d\in\Bbb Z$ with
$2\nmid ac$, $2\mid b$ and $2\mid d$. If $a+bi$ and $c+di$ are
relatively prime elements of $\Bbb Z[i]$, we have the following
general law of quartic reciprocity:
$$\Qs{a+bi}{c+di}=(-1)^{\f b2\cdot \f{c-1}2+\f d2\cdot \f{a+b-1}2}
\Qs{c+di}{a+bi}.$$ In particular, if $4\mid b$, then
$$\Qs{a+bi}{c+di}=(-1)^{\f{a-1}2\cdot\f d2}\Qs{c+di}{a+bi}.$$
\endpro
\pro{Lemma 2.4 ([E], [S1, Lemma 2.1])} Let $a,b,m\in\Bbb Z$ with
$2\nmid m$ and $(m,a^2+b^2)=1$. Then
$$\Qs{a+bi}m^2=\Ls{a^2+b^2}m.$$
\endpro

 \pro{Lemma 2.5 ([S3, Lemma 4.3])} Let $a,b\in\Bbb Z$ with
$2\nmid a$ and $2\mid b$. For any integer $x$ with $(x,a^2+b^2)=1$
we have
$$\Qs{x^2}{a+bi}=\Ls x{a^2+b^2}.$$
\endpro

\par For a given odd prime $p$ let $\Bbb Z_p$ denote the set of
 those rational numbers whose denominator is not divisible by $p$.
 Following [S1,S2] we define
 $$Q_r(p)=\Big\{k\bigm|\ k\in\Bbb Z_p,\
 \Qs{k+i}p=i^r\Big\}\qtq{for}r=0,1,2,3.$$
\pro{Lemma 2.6 ([S1, Theorem 2.3])} Let $p$ be an odd prime,
$r\in\{0,1,2,3\}$, $k\in\Bbb Z_p$ and $k^2+1\not\e 0\mod p$.
\par $(\t{\rm i})$ If $p\e 1\mod 4$ and $t^2\e -1\mod p$ with
 $t\in\Bbb Z_p$, then $k\in Q_r(p)$ if and only if
 $\sls{k+t}{k-t}^{(p-1)/4}\e t^r\mod p$.
\par $(\t{\rm ii})$ If $p\e 3\mod 4$, then $k\in Q_r(p)$ if and only if
 $\sls{k-i}{k+i}^{(p+1)/4}\e i^r\mod p$.\endpro

\pro{Lemma 2.7 ([S5, Lemma 2.9])} Suppose $c,d,m,x\in\Bbb Z$,
$2\nmid m$, $x^2\e c^2+d^2\mod m$ and $(m,x(x+d))=1$. Then
$$\Qs{c+di}m=\Ls{x(x+d)}m.$$
\endpro
\pro{Lemma 2.8} Let $p$ be a prime of the form $8k+1$, $q\in\Bbb Z$,
$2\nmid q$ and $p\nmid q$. Suppose that $p=c^2+d^2=x^2+2qy^2$ with
$c,d,x,y\in\Bbb Z$, $c\e 1\mod 4$, $d=2^rd_0$ and $d_0\e 1\mod 4$.
If $\qs{x/y}{c+di}=(-1)^{s+\f{p-1}8}i^n$, then
$$(-q)^{\f{p-1}8}\e\cases (-1)^s\sls dc^{n-1}\mod p&\t{if $8\mid
d-4$,}
\\(-1)^{s+\f d8}\sls dc^n\mod p&\t{if $8\mid d$.}
\endcases$$
\endpro
Proof. As $c\e 1\mod 4$ and $4\mid d$ we see that $c+di$ is primary
in $\Bbb Z[i]$. Since $i\e d/c\mod {c+di}$ we have
$$\Ls xy^{\f{p-1}4}\e \Qs{x/y}{c+di}=(-1)^{s+\f{p-1}8}i^n
\e (-1)^{s+\f{p-1}8}\Ls dc^n\mod{c+di}.$$ As $\sls xy^2\e -q\mod p$
and the norm of $c+di$ is $p$, from the above we deduce that
$$(-2q)^{\f{p-1}8}\e \Ls xy^{\f{p-1}4}
\e (-1)^{s+\f{p-1}8}\Ls dc^n\mod p.$$ By [L] or [HW, (1.4) and
(1.5)],
$$(-2)^{\f{p-1}8}\e \Ls cd^{-\f d4}
\e \cases \sls cd^{-d_0}\e \sls cd^{-1}=\f dc\mod p&\t{if $8\mid
d-4$,}\\(-1)^{\f d8}\mod p&\t{if $8\mid d$.}
\endcases\tag 2.1$$
Thus the result follows.

 \subheading{3. Congruences for $q^{[p/8]}\mod p$ with
$4p=x^2+qy^2$}
 \pro{Theorem 3.1} Let $p$ be a prime of the form
$4k+1$, $q\in\Bbb Z$, $q\e 3\mod 8$ and $p\nmid q$. Suppose that
$p=c^2+d^2$, $4p=x^2+qy^2$, $c,d,x,y\in\Bbb Z$, $c\e y\e 1\mod 4$,
 $(c,x+2d)=1$ and
$\sqs{2c/(x+2d)+i}q=i^k$. Then
$$(-q)^{[p/8]}\e \cases (-1)^{[\f x4]+\f{q-3}8}\sls dc^{k-1}\mod p&\t{if
$8\mid p-1$,}
\\ -(-1)^{\f{x^2-1}8+\f{q-3}8}\sls dc^{k-1}\f yx\mod p&\t{if $8\mid
p-5$.}\endcases$$
  \endpro

Proof.  Since $4p=x^2+qy^2$ and $y\e 1\mod 4$ we see that $2\nmid
x$, $(x,y)=1$ and $p\nmid x$. Thus $(x,(2c)^2+(x+2d)^2)=(x,4p)=1$.
As $qy^2=(2c)^2+(x+2d)(2d-x)$ we see that $(qy,x+2d)\mid 4c^2$ and
so $(qy,x+2d)=1$. Also,
$$\align &(qy^2,(2c)^2+(x+2d)^2)
\\&=((2c)^2+(x+2d)^2-2x(x+2d),(2c)^2+(x+2d)^2)
=(2x(x+2d),(2c)^2+(x+2d)^2)\\&=(x+2d,(2c)^2+(x+2d)^2)=(x+2d,4c^2)=1.
\endalign$$
Now using Lemmas 2.1-2.5 and the fact that $\sqs an=1$ for
$a,n\in\Bbb Z$ with $2\nmid n$ and $(a,n)=1$ we see that
$$\align i^k&=\Qs{2c+(x+2d)i}q
=\Qs iq\Qs{x+2d-2ci}q\\&=(-1)^{\f{q^2-1}8}\cdot (-1)^{\f{q-1}2} \Qs
q{x+2d-2ci}=\Qs{qy^2}{x+2d-2ci}\Qs{y^2}{x+2d-2ci}
\\&=\Qs{(x+2d)^2+(2c)^2-2x(x+2d)}{x+2d-2ci}\Ls y{(x+2d)^2+4c^2}
\\&=\Qs{-2x(x+2d)}{x+2d-2ci}\Ls y{(x+2d)^2+4c^2}
\\&=-\Qs 2{x+2d-2ci}\Qs{x(x+2d)}{x+2d-2ci}\Ls y{(x+2d)^2+4c^2}
\\&=-i^{(-1)^{(x+1)/2}}(-1)^{\f{x(x+2d)-1}2}\Qs{x+2d-2ci}{x(x+2d)}
\Ls{(x+2d)^2+4c^2}y
\\&=(-1)^{\f{x-1}2}i\Qs{2(d-ci)}x\Qs{-2ci}{x+2d}\Ls{2x(x+2d)+qy^2}y
\\&=(-1)^{\f{x-1}2}i\Qs{d-ci}x\Qs{i}{x+2d}\Ls{2x(x+2d)}y.
\endalign$$
Thus, applying Lemma 2.7 we see that
$$\align i^k&=(-1)^{\f{x-1}2}i\Qs{-i}x\Qs{c+di}x(-1)^{\f{(x+2d)^2-1}8}\cdot
(-1)^{\f{y^2-1}8}\Ls{x(x+2d)}y
\\&=(-1)^{\f{x-1}2}i\cdot (-1)^{\f{x^2-1}8}\Qs{c+di}x
(-1)^{\f{x^2-1}8+\f d2}\cdot(-1)^{\f{4p-x^2-q}8}\Ls{\f x2(\f x2+d)}y
\\&=(-1)^{\f{x-1}2+\f d2}i\cdot (-1)^{\f{x-1}2\cdot \f d2}\Qs
x{c+di} (-1)^{\f{x^2-1}8+\f{q-3}8}\Qs{c+di}y
\\&=-(-1)^{\f{x+1}2(\f d2+1)}i\Qs x{c+di}
(-1)^{\f{x^2-1}8+\f{q-3}8}\Qs y{c+di}
\\&=-(-1)^{\f{x+1}2(\f d2+1)}\cdot(-1)^{\f{x^2-1}8+\f{q-3}8}i\Qs x{c+di}
\Qs{y^{-1}}{c+di}\Qs{y^2}{c+di}
\\&=-(-1)^{\f{x+1}2(\f
d2+1)}\cdot(-1)^{\f{x^2-1}8+\f{q-3}8}i\Qs{x/y}{c+di}\Ls y{c^2+d^2}.
\endalign$$
As $\ls
y{c^2+d^2}=\ls{c^2+d^2}y=\ls{4(c^2+d^2)}y=\ls{x^2+qy^2}y=\ls{x^2}y=1$,
from the above we deduce that
$$\Qs{x/y}{c+di}=-(-1)^{\f{x+1}2(\f
d2+1)}\cdot(-1)^{\f{x^2-1}8+\f{q-3}8}i^{k-1}.$$ Clearly $(-1)^{\f
d2}=(-1)^{\f{p-1}4}$ and $i\e d/c\mod{c+di}$. Since $c+di$ or
$-c-di$ is primary in $\Bbb Z[i]$, we  have
$$\Ls xy^{(p-1)/4}\e
\Qs{x/y}{c+di}\e -(-1)^{\f{x+1}2\cdot\f{p-5}4+\f{x^2-1}8+\f{q-3}8}
\Ls dc^{k-1} \mod{c+di}$$ and so
$$\Ls xy^{(p-1)/4}\e -(-1)^{\f{x+1}2\cdot\f{p-5}4+\f{x^2-1}8+\f{q-3}8}
\Ls dc^{k-1} \mod p.$$ Note that $(x/y)^2\e -q\mod p$. We then have
$$(-q)^{[p/8]}\e \cases \sls xy^{\f{p-1}4}
\e (-1)^{\f{x-1}2+\f{x^2-1}8+\f{q-3}8}\sls dc^{k-1}\mod p&\t{if
$8\mid p-1$,}
\\\sls xy^{\f{p-1}4}\f yx
\e -(-1)^{\f{x^2-1}8+\f{q-3}8}\sls dc^{k-1}\f yx\mod p&\t{if $8\mid
p-5$.}\endcases$$ This yields the result.

\pro{Corollary 3.1} Let $p$ be a prime of the form $12k+1$,
$p=c^2+d^2$, $4p=x^2+27y^2$, $c,d,x,y\in\Bbb Z$, $c\e y\e 1\mod 4$
and $(c,x+2d)=1$.
\par $(\t{\rm i})$ If $p\e 1\mod{24}$, then
$$(-3)^{\f{p-1}8}\e
\cases (-1)^{[\f x4]}\f dc\mod p&\t{if $c\e 0\mod 3$,}
\\\pm (-1)^{[\f x4]}\mod p&\t{if $x-d\e\pm c\mod 3$.}
\endcases$$
\par $(\t{\rm ii})$ If $p\e 13\mod{24}$, then
$$(-3)^{\f{p-5}8}\e
\cases -(-1)^{\f{x^2-1}8}\f {3dy}{cx}\mod p&\t{if $c\e 0\mod 3$,}
\\\pm (-1)^{\f{x^2-1}8}\f{3y}x\mod p&\t{if $x-d\e\pm c\mod 3$.}
\endcases$$
\endpro
Proof. Clearly
$$\Qs i3=-1,\q\Qs{1+i}3=-i,\q \Qs{-1+i}3=i.$$
and
$$(-3)^{[\f p8]} =3^{4[\f p8]}(-27)^{-[\f p8]}\e
3^{(-1)^{(p-1)/4}-1}(-27)^{-[\f p8]}\mod p.$$ Now taking $q=27$ in
Theorem 3.1 and then applying the above we deduce the result.

 \pro{Theorem 3.2} Let
$p$ be a prime of the form $4k+1$, $q\in\Bbb Z$, $q\e 3\mod 8$ and
$p\nmid q$. Suppose that $p=c^2+d^2$, $4p=x^2+qy^2$, $c,d,x,y\in\Bbb
Z$, $c\e y\e 1\mod 4$, $(d,x+2c)=1$ and $\sqs{-2d/(x+2c)+i}q=i^k$.
Then
$$(-q)^{[p/8]}\e \cases (-1)^{[\f x4]} \sls dc^k\mod p&
\t{if $8\mid p-1$,}
\\-(-1)^{\f{x^2-1}8} \sls dc^k\f yx\mod p&\t{if $8\mid
p-5$.}\endcases$$
  \endpro

Proof.  Since $4p=x^2+qy^2$ and $y\e 1\mod 4$ we see that $2\nmid
x$, $(x,y)=1$ and $p\nmid x$. Thus $(x,(2d)^2+(x+2c)^2)=(x,4p)=1$.
As $qy^2=(2d)^2+(x+2c)(2c-x)$ we see that $(qy,x+2c)\mid (2d)^2$ and
so $(qy,x+2c)=1$. Also,
$$\align &(qy^2,(2d)^2+(x+2c)^2)
\\&=((2d)^2+(x+2c)^2-2x(x+2c),(2d)^2+(x+2c)^2)
=(2x(x+2c),(2d)^2+(x+2c)^2)\\&=(x+2c,(2d)^2+(x+2c)^2)=(x+2c,4d^2)=1.
\endalign$$
Now using Lemmas 2.1-2.5 and the fact that $\sqs an=1$ for
$a,n\in\Bbb Z$ with $2\nmid n$ and $(a,n)=1$ we see that
$$\align i^k&=\Qs{-2d+(x+2c)i}q
=\Qs iq\Qs{x+2c+2di}q\\&=(-1)^{\f{q^2-1}8}\Qs
q{x+2c+2di}=-\Qs{qy^2}{x+2c+2di}\Qs{y^2}{x+2c+2di}
\\&=-\Qs{(x+2c)^2+(2d)^2-2x(x+2c)}{x+2c+2di}\Ls y{(x+2c)^2+4d^2}
\\&=-\Qs 2{x+2c+2di}\Qs{x(x+2c)}{x+2c+2di}\Ls y{(x+2c)^2+4d^2}
\\&=-(-1)^{\f d2}\Qs{x+2c+2di}{x(x+2c)}
\Ls{(x+2c)^2+4d^2}y
\\&=-(-1)^{\f{p-1}4}\Qs{2(c+di)}x\Qs{2di}{x+2c}\Ls{2x(x+2c)+qy^2}y
\\&=-(-1)^{\f{p-1}4}\Qs{c+di}x\Qs{i}{x+2c}\Ls{2x(x+2c)}y.
\endalign$$
Thus, applying Lemma 2.7 we see that
$$\align i^k&=-(-1)^{\f{p-1}4}\Qs{c+di}x(-1)^{\f{(x+2c)^2-1}8}
\Ls 2y\Ls{x(x+2c)}y
\\&=-(-1)^{\f{p-1}4}\Qs{c+di}x
(-1)^{\f{x^2-1}8+\f {x+c}2}\Qs iy\Ls{\f x2(\f x2+c)}y
\\&=-(-1)^{\f{p-1}4}\cdot (-1)^{\f{x-1}2\cdot \f d2}\Qs
x{c+di}(-1)^{\f{x^2-1}8+\f {x+1}2} \Qs iy\Qs{d+ci}y
\\&=-(-1)^{(1+\f{x-1}2)(1+\f{p-1}4)}\cdot
(-1)^{\f{x^2-1}8}\Qs x{c+di}\Qs{-c+di}y
\\&=-(-1)^{\f{x+1}2\cdot\f{p+3}4+\f{x^2-1}8}
\Qs x{c+di}\Qs{c+di}y^{-1}
\\&=-(-1)^{\f{x+1}2\cdot\f{p+3}4+\f{x^2-1}8}\Qs x{c+di}\Qs
y{c+di}^{-1}
\\&=-(-1)^{\f{x+1}2\cdot\f{p+3}4+\f{x^2-1}8}\Qs {x/y}{c+di}.\endalign$$
Clearly $(-1)^{\f d2}=(-1)^{\f{p-1}4}$ and $i\e d/c\mod{c+di}$.
Since $c+di$ or $-c-di$ is primary in $\Bbb Z[i]$, we  have
$$\Ls xy^{(p-1)/4}\e
\Qs{x/y}{c+di}\e -(-1)^{\f{x+1}2\cdot\f{p+3}4+\f{x^2-1}8} \Ls dc^k
\mod{c+di}$$ and so
$$\Ls xy^{(p-1)/4}\e -(-1)^{\f{x+1}2\cdot\f{p+3}4+\f{x^2-1}8} \Ls dc^k \mod p.$$
 Note that $(x/y)^2\e -q\mod p$. We then have
$$(-q)^{[p/8]}\e \cases \sls xy^{\f{p-1}4}
\e -(-1)^{\f{x+1}2+\f{x^2-1}8} \sls dc^k=(-1)^{[\f x4]}\sls dc^k\mod
p&\t{if $8\mid p-1$,}
\\\sls xy^{\f{p-1}4}\f yx
\e -(-1)^{\f{x^2-1}8} \sls dc^k\f yx\mod p&\t{if $8\mid
p-5$.}\endcases$$ This is the result.

\par\q
\newline{\bf Remark 3.1}
We note that the $k$ in Theorem 3.1 depends only on $\f c{x+2d}\mod
q$, and  the $k$ in Theorem 3.2 depends only on $\f d{x+2c}\mod q$.

\pro{Theorem 3.3} Let $p$ and $q$ be primes such that $p\e 1\mod 4$
and $q\e 3\mod 8$. Suppose $p=c^2+d^2$, $4p=x^2+qy^2$,
$c,d,x,y\in\Bbb Z$, $c\e y\e 1\mod 4$,  and
$\sls{2(c-di)}x^{\f{q+1}4}\e i^m\mod q$. Assume $(c,x+2d)=1$ or
$(d,x+2c)=1$. Then
$$(-q)^{[p/8]}\e\cases (-1)^{[\f x4]+1}\sls dc^m\mod p&\t{if $p\e 1\mod
8$,}\\(-1)^{\f{x^2-1}8}\sls dc^m\f yx\mod p&\t{if $p\e 5\mod 8$.}
\endcases$$
\endpro
Proof. Clearly $q\nmid x$ and $x$ is odd. We first assume
$(c,x+2d)=1$. By the proof of Theorem 3.1,
$(q,(x+2d)(c^2+(x+2d)^2))=1$. It is easily seen that
$\f{2c/(x+2d)-i}{2c/(x+2d)+i}=\f{2c-(x+2d)i}{2c+(x+2d)i}\e
\f{2(c-di)}{ix}\mod q$.
 Thus, for $k=0,1,2,3$, using Lemma 2.6 we get
$$\align &\Qs{2c/(x+2d)+i}q=i^k\\&\iff\f {2c}{x+2d}\in Q_k(q)
\iff \Ls{\f {2c}{x+2d}-i}{\f {2c}{x+2d}+i}^{\f{q+1}4}\e i^k\mod q
\\&\iff
\Ls{2(c-di)}{ix}^{\f{q+1}4}\e i^k\mod q \iff
\Ls{2(c-di)}x^{\f{q+1}4} \e i^{\f{q+1}4+k}\mod q.
\endalign$$
Since  $\sls{2(c-di)}x^{\f{q+1}4}\e i^m\mod q$, from the above we
deduce that
$$\Qs{2c/(x+2d)+i}q=i^{m-\f{q+1}4}=(-1)^{\f{q-3}8}i^{m-1}.$$ Now,
applying Theorem 3.1 we derive the result.
\par Now we assume $(d,x+2c)=1$. By the proof of Theorem 3.2,
$(q,x+2c)=(q,d^2+(x+2c)^2)=1$.
 It is easily seen
that $\f{2d+(x+2c)i}{2d-(x+2c)i}\e \f{2(c-di)}{-x}\mod q$.
 Thus, for $k=0,1,2,3$, using Lemma 2.6 we get
$$\align &\Qs{-2d/(x+2c)+i}q=i^k\\&\iff-\f {2d}{x+2c}\in Q_k(q)
\iff \Ls{-\f {2d}{x+2c}-i}{-\f {2d}{x+2c}+i}^{\f{q+1}4}\e i^k\mod q
\\&\iff \Ls{2d+(x+2c)i}{2d-(x+2c)i}^{\f{q+1}4}\e i^k\mod q\iff
\Ls{2(c-di)}{-x}^{\f{q+1}4}\e i^k\mod q
\\&\iff \Ls{2(c-di)}x^{\f{q+1}4}\e i^{\f{q+1}2+k}\mod q.
\endalign$$
Since  $\sls{2(c-di)}x^{\f{q+1}4}\e i^m\mod q$, from the above we
deduce that
$\sqs{-2d/(x+2c)+i}q=i^{m-\f{q+1}2}=(-1)^{\f{q+1}4}i^m=-i^m$. Now
applying Theorem 3.2 we deduce the result. The proof is now
complete.

\pro{Theorem 3.4} Let $p$ and $q$ be primes such that $p\e 1\mod 4$
and $q\e 3\mod 8$. Suppose $p=c^2+d^2$, $4p=x^2+qy^2$,
$c,d,x,y\in\Bbb Z$, $c\e y\e 1\mod 4$,  and $q\mid cd$. Assume that
$(c,x+2d)=1$ or $(d,x+2c)=1$. Then
$$(-q)^{[\f p8]}\e\cases \pm(-1)^{[\f x4]+1}\mod p&\t{if $p\e 1\mod
8$ and $x\e \pm 2c\mod q$,}\\\pm (-1)^{\f{q-3}8+[\f x4]}\f dc\mod
p&\t{if $p\e 1\mod 8$ and $x\e \pm 2d\mod q$,}
\\\pm(-1)^{\f{x^2-1}8}\f yx\mod p&\t{if $p\e 5\mod
8$ and $x\e \pm 2c\mod q$,}\\\mp (-1)^{\f{q-3}8+\f{x^2-1}8} \f
{dy}{cx}\mod p&\t{if $p\e 5\mod 8$ and $x\e \pm 2d\mod q$,}
\endcases$$
\endpro
Proof. As $\sls x2^2\e c^2+d^2\mod q$ we see that $x\e \pm 2c\mod q
\Leftrightarrow q\mid d$ and that $x\e \pm 2d\mod q \Leftrightarrow
q\mid c$. Thus,
$$\Qs{c-di}{x/2}^{\f{q+1}4}
\e\cases (\pm 1)^{\f{q+1}4}=\pm 1\mod q&\t{if $x\e \pm 2c\mod q$,}
\\(\mp i)^{\f{q-1}4}=\mp (-1)^{\f{q-3}8}i\mod q&\t{if $x\e \pm
2d\mod q$.}
\endcases$$
Now applying Theorem 3.3 we derive the result.

\pro{Theorem 3.5} Let $p$ be a prime of the form $4k+1$ and
$\sls{-11}p=1$, $p=c^2+d^2$, $4p=x^2+11y^2$, $c,d,x,y\in\Bbb Z$,
$c\e 1\mod 4$, $d=2^rd_0(2\nmid d_0)$, $y=2^ty_0$ and $y_0\e 1\mod
4$. Assume that $(c,x+2d)=1$ or $(d_0,x+2c)=1$.
\par $(\t{\rm i})$ If $p\e 1\mod 8$, then
$$(-11)^{[\f p8]}\e
\cases \pm(-1)^{[\f x4]}\mod p&\t{if $2\nmid x$ and $x\e \pm 4c,\pm
9c\mod{11}$,}
\\\pm (-1)^{[\f x4]}\f dc\mod p&\t{if $2\nmid x$
 and $x\e \pm 4d,\pm 9d\mod{11}$,}
\\\mp(-1)^{[\f x8]+\f y8}\mod p&\t{if $2\mid x$ and $x\e \pm 4c,\pm
9c\mod{11}$,}
\\\mp (-1)^{[\f x8]+\f y8}\f dc\mod p&\t{if $2\mid x$ and $x\e \pm 4d,\pm 9d\mod{11}$.}
\endcases$$

\par $(\t{\rm ii})$ If $p\e 5\mod 8$, then
$$(-11)^{[\f p8]}\e
\cases \mp (-1)^{\f{x^2-1}8}\f yx\mod p&\t{if  $2\nmid x$ and $x\e
\pm 4c,\pm 9c\mod{11}$,}
\\\mp (-1)^{\f {x^2-1}8}\f {dy}{cx}\mod p&\t{if $2\nmid x$ and
$x\e \pm 4d,\pm 9d\mod{11}$,}
\\\mp(-1)^{\f {p-5}8}\f yx\mod p&\t{if $2\mid x$ and
$x\e \pm 4c,\pm
9c\mod{11}$,}
\\\mp (-1)^{\f {p-5}8}\f{dy}{cx}\mod p&\t{if $2\mid x$
and $x\e \pm 4d,\pm 9d\mod{11}$.}
\endcases$$
\endpro
Proof. As $$\ls x2^2\e c^2+d^2\mod{11}\qtq{and}
(c-di)^3=c(c^2-3d^2)+d(d^2-3c^2)i,$$ we see that
$$\Qs{c-di}{x/2}^3\e\cases \mp 1\mod p&\t{if $x\e \pm 4c,\pm
9c\mod{11}$,}
\\\mp i\mod p&\t{if $x\e \pm 4d,\pm
9d\mod{11}$.}\endcases$$ When $2\nmid x$, from the above and Theorem
3.3 (with $q=11$) we deduce the result. When $2\mid x$ and $p\e
1\mod 8$, we have $8\mid y$ and so $(-1)^{\f{p-1}8+\f{x/2-1}2}
=(-1)^{\f{(x/2)^2-1}8+\f{x/2-1}2}=(-1)^{[\f x8]}$. Thus, applying
the above and [S5, Theorem 4.1 (with $q=11$)] we obtain the result.

\par\q
\newline{\bf Remark 3.2}  Let $p$ be a prime of the form $4k+1$,
$q\in\Bbb Z$, $q\e 3\mod 8$, $p\nmid q$,  $p=c^2+d^2$,
$4p=x^2+qy^2$, $c,d,x,y\in\Bbb Z$, $2\nmid xy$, $c\e  1\mod 4$, we
conjecture that one can always choose the sign of $x$ such that
$(c,x+2d)=1$ or $(d,x+2c)=1$. Thus the condition $(c,x+2d)=1$ or
$(d,x+2c)=1$ in Theorems 3.1-3.5 can be canceled. See also [S5,
Remark 4.1].

\subheading{4. Congruences for $(-q)^{(p-1)/8}\mod p$ with
$p=x^2+2qy^2$}
 \pro{Theorem 4.1} Let $p$ be a prime of the form
$8k+1$, $q\in\Bbb Z$, $2\nmid q$ and $p\nmid q$. Suppose that
$p=c^2+d^2=x^2+2qy^2$ with $c,d,x,y\in\Bbb Z$, $c\e 1\mod 4$,
$d=2^rd_0$, $d_0\e 1\mod 4$ and $(c,x+d)=1$. Assume that
$\sqs{c/(x+d)+i}q=i^k$. Then
$$(-q)^{\f{p-1}8}\e
\cases (-1)^{\f {q-1}8+\f d4+\f y2} \sls dc^k\mod p&\t{if $q\e 1\mod
8$,}
\\(-1)^{\f {q-3}8+\f{x-1}2}\sls dc^{k+1}\mod p
 &\t{if $q\e  3\mod 8$,}
\\(-1)^{\f {q-5}8+ \f d4+\f{x-1}2+\f
y2}\sls dc^{k-1}\mod p
 &\t{if $q\e  5\mod 8$,}
\\ (-1)^{\f {q+1}8} \sls dc^k\mod p&\t{if $q\e 7\mod 8$.}
 \endcases$$
\endpro
Proof. Suppose $y=2^ty_0$ $(2\nmid y_0)$. We may assume $y_0\e 1\mod
4$. Since $p=c^2+d^2=x^2+2qy^2\e 1\mod 8$ we see that $2\nmid x$,
$2\mid y$, $4\mid d$, $(x,y)=1$ and $p\nmid x$. Thus
$(x,c^2+(x+d)^2)=(x,p)=1$. As $2qy^2=c^2+(d+x)(d-x)$ we see that
$(qy,x+d)\mid c^2$ and so $(qy,x+d)=1$. Also,
$$\align &(qy^2,(c^2+(x+d)^2)/2)
\\&=((c^2+(x+d)^2)/2-x(x+d),(c^2+(x+d)^2)/2)
=(x(x+d),(c^2+(x+d)^2)/2)\\&=(x(x+d),c^2+(x+d)^2)=(x+d,c^2)=1.
\endalign$$
 It is
easily seen that
$$c+(x+d)i=i^{\f{1\mp 1}2}(1+i)\Big(\f{x+d\pm
c}2+\f{\pm(x+d)-c}2i\Big)$$ and so $$\Ls{x+d\pm
c}2^2+\Ls{\pm(x+d)-c}2^2=\f {c^2+(x+d)^2}2.$$ Set
$\ep=(-1)^{\f{x-1}2}$. As $4\mid d$ and $4\mid c-1$ we have $x+d\e
\ep\mod 4$ and $4\mid (\ep(x+d)-c)$. Using Lemmas 2.1-2.5 and [S5,
Lemma 2.10(ii)] and the above we see that
$$\align i^k&
=\Qs{c+(x+d)i}q=\Qs iq^{\f{1-\ep}2}\Qs{1+i}q\Qs{\f{x+d+\ep
c}2+\f{\ep(x+d)-c}2i}q
\\&=(-1)^{\f{q-\sls{-1}q}4\cdot\f{x-1}2}i^{\f{\sls{-1}qq-1}4}
(-1)^{\f{q-1}2\cdot\f{\ep(x+d)-c}4}\Qs q{\f{x+d+\ep
c}2+\f{\ep(x+d)-c}2i}
\endalign$$
and
$$\align &\Qs q{\f{x+d+\ep
c}2+\f{\ep(x+d)-c}2i}
\\&=\Qs {qy^2}{\f{x+d+\ep
c}2+\f{\ep(x+d)-c}2i}\Qs {y^2}{\f{x+d+\ep c}2+\f{\ep(x+d)-c}2i}
\\&=\Qs {(c^2+(x+d)^2)/2-x(x+d)}{\f{x+d+\ep
c}2+\f{\ep(x+d)-c}2i}\Ls y{\sls{x+d+\ep c}2^2+\sls{\ep(x+d)-c}2^2}
\\&=\Qs {-x(x+d)}{\f{x+d+\ep
c}2+\f{\ep(x+d)-c}2i}\Ls y{(c^2+(x+d)^2)/2}
\\&=(-1)^{\f{\ep(x+d)-c}4}\Qs {\f{x+d+\ep
c}2+\f{\ep(x+d)-c}2i}{x(x+d)}
 (-1)^{\f{c^2-(x+d)^2}8t+\f d4t}\Qs{y^{-1}}{c+di}.
\endalign$$
Clearly, $$(-1)^{\f{c^2-(x+d)^2}8}=(-1)^{\f{c^2-x^2-2dx}8}
=(-1)^{\f{c-\ep x}4\cdot \f{c+\ep x}2-\f d4\ep}=(-1)^{\f{c-\ep
x}4-\f d4\ep}=(-1)^{\f{\ep(x+d)-c}4}.$$ It is easily seen that
$$\align&\Qs {\f{x+d+\ep
c}2+\f{\ep(x+d)-c}2i}{x(x+d)}
\\&=\Qs {x+d+\ep
c+(\ep(x+d)-c)i}{x}\Qs {x+d+\ep c+(\ep(x+d)-c)i}{x+d}
\\&=\Qs{d+\ep c+(\ep d-c)i}x\Qs{\ep c-ci}{x+d}=\Qs{(\ep-i)(c+di)}x
\Qs{\ep -i}{x+d}
\\&=\Qs{\ep -i}{x(x+d)}\Qs{c+di}x=
\Qs{i^{\f{5+\ep}2}(1+i)}{x(x+d)}\Qs{c+di}x
\\&=\Qs i{x(x+d)}^{\f{5+\ep}2}\Qs{1+i}{x(x+d)}\Qs x{c+di}
\\&=(-1)^{\f{x(x+d)-1}4\cdot\f{5+\ep}2}
i^{\f{x(x+d)-1}4} \Qs x{c+di} \\&=(-1)^{\f
d4\cdot\f{x+1}2+\f{x^2-1}8}i^{\f{dx}4}\Qs x{c+di}.\endalign$$
Hence
$$\align\Qs q{\f{x+d+\ep c}2+\f{\ep(x+d)-c}2i}
&=(-1)^{\f{\ep(x+d)-c}4(1+t)+\f d4t}\cdot (-1)^{\f
d4\cdot\f{x+1}2+\f{x^2-1}8}i^{\f{dx}4}\Qs {x/y}{c+di}
\\&=(-1)^{\f{\ep x-c}4(1+t)+\f d4\cdot\f{x-1}2+\f{x^2-1}8}
i^{\f{dx}4}\Qs {x/y}{c+di}.\endalign$$
 Therefore
$$\align i^k
&=(-1)^{\f{q-\sls{-1}q}4\cdot\f{x-1}2 +\f{q-1}2(\f{\ep x-c}4+\f d4)}
i^{\f{\sls{-1}qq-1}4}
\\&\q\times (-1)^{\f{\ep x-c}4(1+t)+\f d4\cdot\f{x-1}2+\f{x^2-1}8}
i^{\f{dx}4}\Qs {x/y}{c+di}.\endalign$$
It is clear that
$$ \align&i^{\f {dx}4}=i^{\f d4(x-1)+\f d4}=(-1)^{\f
d4\cdot\f{x-1}2}i^{\f d4},\
(-1)^{\f{x^2-1}8}=(-1)^{\f{p-1-2qy^2}8}=(-1)^{\f{p-1}8+\f y2},
\\&(-1)^{\f{\ep x-c}4}=(-1)^{\f{\ep x-c}4\cdot
 \f{\ep x+c}2}
=(-1)^{\f{x^2-c^2}8}=(-1)^{\f{d^2-8q(\f y2)^2}8}=(-1)^{\f
y2}\endalign$$ and so $(-1)^{\f{\ep x-c}4(1+t)}=(-1)^{\f
y2(1+t)}=1$. Thus,
$$\Qs{x/y}{c+di}=(-1)^{\f{q-\sls{-1}q}4\cdot\f{x-1}2+\f{q-1}2(\f
y2+\f d4)+\f{p-1}8+\f y2}i^{k-\f d4-\f{q\sls{-1}q-1}4}.$$ Now
applying Lemma 2.8 we deduce the result.

\pro{Theorem 4.2} Let $p$ be a prime of the form $4k+1$, $q\in\Bbb
Z$, $2\nmid q$ and $p\nmid q$. Suppose that $p=c^2+d^2=x^2+2qy^2$
with $c,d,x,y\in\Bbb Z$, $c\e 1\mod 4$, $d=2^rd_0$, $d_0\e 1\mod 4$,
$(d_0,x+c)=1$ and $\sqs{-d/(x+c)+i}q=i^k$. Then
$$(-q)^{\f{p-1}8}\e
\cases (-1)^{\f{x+1}2\cdot\f{q-1}4+\f d4+\f y2}\sls dc^k\mod p&\t{if
$q\e 1\mod 4$,}
\\(-1)^{\f{x+1}2\cdot\f{q+1}4}\sls dc^k\mod p
&\t{if $q\e 3\mod 4$.}
\endcases$$
  \endpro
  Proof. Suppose  $2^m\ \Vert\ (x+c)$ and $y=2^ty_0$ $(2\nmid y_0)$.
   We may assume $y_0\e 1\mod 4$. As $(x+c,d_0)=1$, by Lemma 2.11 we have $(qy_0,x+c)=1$
and $(qy_0^2,(x+c)^2+d^2)=1$. Note that $(x,y)^2\mid p$. We also
have $(x,y)=1$. We prove the theorem by considering the following
three cases: \par {\bf Case 1.} $m<r$. Using Lemmas 2.1-2.5 and the
fact $\qs aq=1$ for $a\in\Bbb Z$ with $(a,q)=1$ we see that
$$\aligned \Qs{d-(x+c)i}q&=\Qs {-2^m}q\Qs iq\Qs{\f{x+c}{2^m}+\f d{2^m}i}q
=(-1)^{\f{q^2-1}8+\f{q-1}2\cdot\f d{2^{m+1}}}\Qs q{\f{x+c}{2^m}+\f
d{2^m}i}
\\&=(-1)^{\f{q^2-1}8+\f{q-1}2\cdot
\f d{2^{m+1}}}\Qs{qy^2}{\f{x+c}{2^m}+\f d{2^m}i}\Qs
{y^2}{\f{x+c}{2^m}+\f d{2^m}i}
\\&=(-1)^{\f{q^2-1}8+\f{q-1}2\cdot
\f d{2^{m+1}}}\Qs{((x+c)^2+d^2)/2-x(x+c)}{\f{x+c}{2^m}+\f d{2^m}i}
\Ls y{\f{(x+c)^2+d^2}{2^{2m}}}
\\&=(-1)^{\f{q^2-1}8+\f{q-1}2\cdot
\f d{2^{m+1}}}\Qs{-2^mx(x+c)/2^m}{\f{x+c}{2^m}+\f d{2^m}i} \Ls
y{\f{(x+c)^2+d^2}{2^{2m}}} .\endaligned$$ By [S5, p.15],
$$\Ls{y}{((x+c)^2+d^2)/2^{2m}}
=(-1)^{\f d{2^{m+1}}t+\f d4t}\Qs{y^{-1}}{c+di}.$$ We also have
$$\Qs {2^m}{\f{x+c}{2^m}+\f d{2^m}i}=\Qs 2{\f{x+c}{2^m}+\f
d{2^m}i}^m =i^{\f{x+c}{2^m}\cdot\f {dm}{2^{m+1}}}$$ and
$$\align \Qs{-x(x+c)/2^m}{\f{x+c}{2^m}+\f d{2^m}i}
&=(-1)^{\f{-x(x+c)/2^m-1}2\cdot\f d{2^{m+1}}} \Qs {\f{x+c}{2^m}+\f
d{2^m}i}{-x(x+c)/2^m}
\\&=(-1)^{\f{x(x+c)/2^m+1}2\cdot\f d{2^{m+1}}}
 \Qs {\f{x+c}{2^m}+\f
d{2^m}i}{-x}\Qs  {\f{x+c}{2^m}+\f d{2^m}i}{(x+c)/2^m}
\\&=(-1)^{\f{x(x+c)/2^m+1}2\cdot\f d{2^{m+1}}}
\Qs{x+c+di}x\Qs{\f d{2^m}i}{(x+c)/2^m}
\\&=(-1)^{\f{x(x+c)/2^m+1}2\cdot\f d{2^{m+1}}}
\Qs{c+di}x\Qs{i}{(x+c)/2^m}
\\&=(-1)^{\f{x(x+c)/2^m+1}2\cdot\f d{2^{m+1}}}\Qs
x{c+di}(-1)^{\f{\sls{x+c}{2^m}^2-1}8}.\endalign$$ Thus,
$$\aligned i^k&=\Qs{d-(x+c)i}q=(-1)^{\f{q^2-1}8+\f{q-1}2\cdot\f
d{2^{m+1}}} i^{\f{x+c}{2^m}\cdot\f {dm}{2^{m+1}}}
\\&\q\times(-1)^{\f{x(x+c)/2^m+1}2\cdot
\f d{2^{m+1}}+\f{\sls{x+c}{2^m}^2-1}8} (-1)^{\f d{2^{m+1}}t+\f
d4t}\Qs{x/y}{c+di}.\endaligned\tag 4.1$$ As
$2qy^2=d^2-(x+c)^2+2c(x+c)$ we have
$$q\f{y^2}{2^m}=2^{2r-m-1}d_0^2-2^{m-1}\Ls{x+c}{2^m}^2+c\cdot
\f{x+c}{2^m}.\tag 4.2$$
\par Suppose $x\e 1\mod 4$. Then $m=1<r$. From (4.2) we see that
$ 2^{2t-1}q\e 2^{2r-2}-1+c\cdot\f{x+c}2\mod 8$ and so $\f{x+c}2\e
c(2^{2t-1}q-2^{2r-2}+1)\mod 8.$ If $8\nmid d$, then $r=2$ and
$\f{x+c}2\e c(2^{2t-1}q-3)\mod 8$. Thus,
$$\align (-1)^{\f{(\f{x+c}2)^2-1}8}
&=(-1)^{\f{c^2-1}8+\f{(2^{2t-1}q-3)^2-1}8}
 =(-1)^{\f{c^2-1}8+\f{(2^{2t-2}q-2)(2^{2t-2}q-1)}2}
 \\&=(-1)^{\f{c^2-1}8+1+\f{q+1}2\cdot \f y2}=(-1)^{\f{p-1}8+1+\f{q+1}2\cdot \f y2}\endalign$$
and
$$(-1)^{\f{(x+c)/2+1}2}i^{\f{x+c}2}=(-1)^{2^{2t-2}q-1}i^{2^{2t-1}q-3}
=-(-1)^{2^{2t-2}}\cdot (-1)^{2^{2t-2}}i=-i.$$ Hence, from (4.1) we
deduce that
$$\align i^k&=(-1)^{\f{q^2-1}8+\f{q-1}2}i^{\f{x+c}2}(-1)^{\f{(x+c)/2+1}2}
\cdot (-1)^{\f{\sls{x+c}2^2-1}8}\Qs{x/y}{c+di}
\\&=(-1)^{\f{q^2-1}8+\f{q-1}2}(-i)(-1)^{\f{p-1}8+1+\f{q+1}2
\cdot \f y2}\Qs{x/y}{c+di}.
\endalign$$
That is,
$$\Qs{x/y}{c+di}
=(-1)^{\f{q^2-1}8+\f{q-1}2+\f{p-1}8+\f{q+1}2\cdot\f y2}i^{k-1}.$$
Now applying Lemma 2.8 we obtain the result.
\par If $8\mid d$, from (4.2) we see that $2^{2t-1}q\e 2q(\f y2)^2\e
-1+c\cdot \f{x+c}2\mod 8$ and so $\f{x+c}2\e c(2^{2t-1}q+1)\mod 8$.
Thus,
$$\align (-1)^{\f{(\f{x+c}2)^2-1}8}&=(-1)^{\f{c^2-1}8+\f{(2^{2t-1}q+1)^2-1}8}
=(-1)^{\f{c^2+d^2-1}8+\f{2^{2t-2}q(2^{2t-2}q+1)}2}
\\&=(-1)^{\f{p-1}8 +\f{q+1}2\cdot \f y2}.\endalign$$
Now, from (4.1) and the above we derive that
$$\align i^k&=(-1)^{\f{q^2-1}8}i^{\f{x+c}2\cdot \f d4}
(-1)^{\f{(\f{x+c}2)^2-1}8}\Qs{x/y}{c+di}
\\&=(-1)^{\f{q^2-1}8+\f d8+\f{p-1}8+\f{q+1}2\cdot \f
y2}\Qs{x/y}{c+di}.\endalign$$ Now applying Lemma 2.8 we obtain
$$(-q)^{\f{p-1}8}\e (-1)^{\f{q^2-1}8+\f{q+1}2\cdot\f y2}\Ls dc^k
=\cases
 (-1)^{\f{q-1}4+\f d4+\f y2}\sls dc^k\mod p&\t{if $4\mid q-1$,}
 \\\sls dc^k\mod p&\t{if $4\mid q-3$.}\endcases$$
This is the result.
\par
Now assume $r>m\ge 2$. Then $x\e 3\mod 4$, $2r-m-1\ge
2(m+1)-m-1=m+1\ge 3$ and so $q\f{y^2}{2^m}\e
-2^{m-1}+c\cdot\f{x+c}{2^m}\mod 8$. Hence $2^m\ \Vert\ y^2$, $m=2t$
and so $q\e -2^{m-1}+c\cdot\f{x+c}{2^m}\mod 8$. That is,
$\f{x+c}{2^m}\e c(2^{m-1}+q)\mod 8$. Thus,
$$\align &(-1)^{\f{q^2-1}8+\f{q-1}2\f
d{2^{m+1}}}(-1)^{\f{x(x+c)/2^m+1}2\cdot \f
d{2^{m+1}}+\f{\sls{x+c}{2^m}^2-1}8}\cdot (-1)^{\f d{2^{m+1}}t+\f
d4t} i^{\f{x+c}{2^m}\cdot\f {dm}{2^{m+1}}}
\\&=(-1)^{\f{q^2-1}8+\f{q-1}2\cdot\f
d{2^{m+1}}}(-1)^{\f {-(2^{m-1}+q)+1}2\cdot
2^{r-m-1}d_0+\f{c^2(2^{m-1}+q)^2-1}8}(-1)^{\f d{2^{m+1}}t+\f d4t}
(-1)^{\f{dt}{2^{m+1}}}
\\&=(-1)^{\f{q^2-1}8+\f{q-1}2\cdot 2^{r-m-1}d_0}\cdot(-1)^{(2^{m-2}+\f
{q-1}2)2^{r-m-1}+\f{c^2-1}8+\f{(2^{m-1}+q)^2-1}8}\cdot (-1)^{\f d4t}
\\&=(-1)^{\f{q^2-1}8+\f{q-1}2\cdot 2^{r-m-1}}\cdot
(-1)^{(2^{m-2}+\f{q-1}2)2^{r-m-1}+\f{p-1}8+\f{q^2-1}8+2^{m-3}
(2^{m-2}+q)}
\\&=(-1)^{2^{r-3}+\f{p-1}8+2^{m-3}
(2^{m-2}+q)} =\cases (-1)^{\f d8+\f{p-1}8+\f{q+1}2}&\t{if $m=2$,}
\\(-1)^{\f{p-1}8}&\t{if $m>2$.}\endcases
\endalign$$
Hence, from (4.1) and the above we get
$$\Qs{x/y}{c+di}=
\cases (-1)^{\f d8+\f{p-1}8+\f{q+1}2}i^k&\t{if $r>m=2$,}
\\(-1)^{\f{p-1}8}i^k&\t{if $r>m>2$.}\endcases$$
Now applying Lemma 2.8 and the fact $m=2t$ we obtain
 $$\aligned(-q)^{\f{p-1}8}&\e
\cases (-1)^{\f{q+1}2}\sls dc^k\mod p&\t{if $r>m=2$,}
\\(-1)^{\f d8}\sls dc^k=\sls dc^k\mod p&\t{if $r>m>2$}
\endcases
\\&=(-1)^{\f{q+1}2\cdot\f y2}\Ls dc^k
=(-1)^{\f{q+1}2(\f d4+\f y2)}\Ls dc^k \mod p\qtq{for}r>m\ge 2.
\endaligned$$ This yields the result.
 \par {\bf Case 2.} $m=r$. As $p\e 1\mod 8$ we have $4\mid d$,
$m=r\ge 2$, $4\mid x+c$ and so $x\e 3\mod 4$.
 It is clear that
$$\aligned i^k&=\Qs{-d/(x+c)+i}q=\Qs{-\f d{2^r}+\f{x+c}{2^r}i}q
\\&=\Qs{1+i}q\Qs {\f{x+c-d}{2^{r+1}}+\f{x+c+d}{2^{r+1}}i}q
=\Qs{1+i}q\Qs {\f{x+c-d}{2^{r+1}}-\f{x+c+d}{2^{r+1}}i}q^{-1} \\&=\Qs
iq\Qs{\f{x+c}{2^r} +\f d{2^r}i}q =\Ls 2q\Qs{1+i}q\Qs
{\f{x+c+d}{2^{r+1}}-\f{x+c-d}{2^{r+1}}i}q.\endaligned\tag 4.3$$ As
$2^r\ \Vert\ x+c$ and $2^r\ \Vert\ d$ we have $2^{r+1}\mid x+c\pm d$
and $\f{d+x+c}{2^r}\cdot\f{d-x-c}{2^r}=\sls
d{2^r}^2-\sls{x+c}{2^r}^2 \e 0\mod 8$. Since
$2qy^2=(d+x+c)(d-x-c)+2c(x+c)$ and $r\ge 2$ we see that
$$c\f{x+c}{2^r}=\f{qy^2}{2^r}-\f{(d+x+c)(d-x-c)}{2^{r+1}}\e
q\f{y^2}{2^r}\mod {2^{r+2}}.\tag 4.4$$ Thus, $r=2t$, $(-1)^{\f d4+\f
y2}=1$ and $\f{x+c}{2^r}\e cq\mod 8$. Hence
$$\f{x+c}{2^r}\e q\e (-1)^{\f{q-1}2}\e (-1)^{\f{q-1}2}d_0\e
(-1)^{\f{q-1}2}\f d{2^r}\mod 4.$$ Set
$$A=\f{x+c+(-1)^{\f{q-1}2}d}{2^{r+1}}\qtq{and}
B=\f{x+c-(-1)^{\f{q-1}2}d}{2^{r+1}}.$$ Then
$$2\nmid A,\ 2\mid B\qtq{and}A^2+B^2=\f{(x+c)^2+d^2}{2^{2r+1}}.$$
Thus,
$$\align \Qs{A-Bi}q&=(-1)^{\f{q-1}2\cdot\f B2}\Qs q{A-Bi}
=(-1)^{\f{q-1}2\cdot\f B2}\Qs{qy^2}{A-Bi}\Qs{y^2}{A-Bi}
\\&=(-1)^{\f{q-1}2\cdot\f B2}\Qs{((x+c)^2+d^2)/2-x(x+c)}{A-Bi}
\Ls y{A^2+B^2}
\\&=(-1)^{\f{q-1}2\cdot\f B2}\Qs{-2^rx(x+c)/2^r}{A-Bi}\Ls
y{((x+c)^2+d^2)/2^{2r+1}}
\\&=(-1)^{\f{q-1}2\cdot\f B2}\Qs 2{A-Bi}^r\Qs{-x(x+c)/2^r}{A-Bi}
\Ls {2^ty_0}{((x+c)^2+d^2)/2^{2r+1}}
\\&=(-1)^{\f{q-1}2\cdot\f
B2}i^{-\f{AB}2r}(-1)^{\f{-x(x+c)/2^r-1}2\cdot \f B2}
\Qs{A-Bi}{-x(x+c)/2^r} \\&\q\times(-1)^{\f{\f
12(\sls{x+c}{2^r}^2+\sls
d{2^r}^2)-1}4t}\Ls{((x+c)^2+d^2)/2^{2r+1}}{y_0} .\endalign$$ Since
$x\e 3\mod 4$ and $\f{x+c}{2^r}\e q\mod 4$ we have
$(-1)^{\f{-x(x+c)/2^r-1}2}=(-1)^{\f{q-1}2}$. As $\f{x+c}{2^r}\e
cq\mod 8$, we see that $\f B2=\f 14(\f{x+c}{2^r}-(-1)^{\f{q-1}2}\f
d{2^r}) \e \f {cq-(-1)^{\f{q-1}2}d_0}4\mod 2$ and
$$\align (-1)^{\f{\f
12(\sls{x+c}{2^r}^2+\sls
d{2^r}^2)-1}4}&=(-1)^{\f{\sls{x+c}{2^r}^2-1}8 +\f{\sls d{2^r}^2-1}8}
=(-1)^{\f{(cq)^2-1}8+\f{d_0^2-1}8} \\&=(-1)^{\f{c^2-1}8+\f{q^2-1}8
+\f{d_0^2-1}8}=(-1)^{\f{p-1}8+\f{q^2-1}8 +\f{d_0^2-1}8}.
\endalign$$
Recall that $r=2t$. From the above we deduce that
$$\align &\Qs{A-Bi}q\\&=(-1)^{\f{AB}2t}\cdot
(-1)^{(\f{p-1}8+\f{q^2-1}8+\f{d_0^2-1}8)t} \Qs{A-Bi}{x(x+c)/2^r}
\Ls{((x+c)^2+d^2)/2}{y_0}
\\&=(-1)^{(\f {cq-(-1)^{(q-1)/2}d_0}4+\f{p-1}8+\f{q^2-1}8+\f{d_0-1}4)t}
\Qs{A-Bi}{x(x+c)/2^r} \Ls{((x+c)^2+d^2)/2}{y_0}
\\&=(-1)^{(\f {cq-(-1)^{(q-1)/2}}4+\f{p-1}8+\f{q^2-1}8)t}
\\&\q\times \Qs{x+c+(-1)^{\f{q-1}2}d-(x+c-(-1)^{\f{q-1}2}d)i}
{x(x+c)/2^r}\Ls{qy^2+x(x+c)}{y_0}
\\&=(-1)^{(\f{c-1}4+\f{q-(-1)^{(q-1)/2}}4+\f{p-1}8+\f{q^2-1}8)t}
\\&\q\times \Qs{c+(-1)^{\f{q-1}2}d-(c-(-1)^{\f{q-1}2}d)i}x
\Qs{(-1)^{\f{q-1}2}d+(-1)^{\f{q-1}2}di} {(x+c)/2^r} \Ls{x(x+c)}{y_0}
\\&=\Qs{-i(1+i)(c+(-1)^{\f{q-1}2}di)}x
\Qs {1+i}{(x+c)/2^r} \Ls{x(x+c)}{y_0} .\endalign$$ Recall that
$\f{x+c}{2^r}\e cq\mod 8$. From the above and  Lemma 2.7 we deduce
that
$$\align \Qs{A-Bi}q
&=\Qs ix\Qs{1+i}x\Qs{c+(-1)^{\f{q-1}2}di}x
i^{\f{(-1)^{\f{(x+c)/2^r-1}2}\f{x+c}{2^r}-1}4} \Ls{x(x+c)}{y_0}
\\&=(-1)^{\f{x+1}4}\cdot i^{\f{-x-1}4}
\Qs
x{c+(-1)^{\f{q-1}2}di}i^{\f{(-1)^{(q-1)/2}\f{x+c}{2^r}-1}4}\Qs{d+ci}{y_0}
\\&\endalign$$
Using (4.4) we see that
$$\align i^{\f{(-1)^{\f{q-1}2}\f{x+c}{2^r}-1}4}
&=i^{\f{(-1)^{\f{q-1}2}c\f{x+c}{2^r}-c}4}
=i^{\f{(-1)^{\f{q-1}2}qy_0^2-c}4}
=i^{\f{\sls{-1}qq(y_0^2-1)+\sls{-1}q-c}4}
\\&=(-1)^{\f{y_0^2-1}8}i^{\f{\sls {-1}qq-c}4}
=\Qs i{y_0}i^{\f{\sls {-1}qq-c}4}.\endalign$$ Also,
 $$i^{-\f{x+1}4}=i^{-\f{x+c}4}\cdot i^{\f{c-1}4}
 =\cases i^{\f{c-1}4}&\t{if $r>2$,}
 \\i^{-q+\f{c-1}4}=(-1)^{\f{q+1}2}i^{1+\f{c-1}4}&\t{if $r=2$}
 \endcases$$
 and
$$\align \Qs i{y_0}\Qs{d+ci}{y_0}&
=\Qs{c-di}{y_0}=\Qs{y_0}{c-di} =\Qs 2{c-di}^{-t}\Qs y{c+di}
\\&=i^{\f d2t}\Qs y{c-di}
=(-1)^{\f d4\cdot \f r2}\Qs y{c-di}
\\&=(-1)^{\f d4}\Qs y{c+di}^{-1}
=(-1)^{\f d4}\Qs {y^{-1}}{c+di}.\endalign$$ Thus,
$$\align \Qs{A-Bi}q&=(-1)^{\f{x+1}4}i^{-\f{x+1}4+\f{\sls{-1}qq-c}4}
\Qs x{c+(-1)^{\f{q-1}2}di}(-1)^{\f d4}\Qs {y^{-1}}{c+di}
\\&=(-1)^{\f{x+1}4+\f d4}i^{-\f{x+c}4+\f{\sls{-1}qq-1}4}
\Qs x{c+(-1)^{\f{q-1}2}di}\Qs {y^{-1}}{c+di}.
\endalign$$
For $r=2$ we have  $\f{x+c}4\e q\mod 4$ and so $i^{-\f{x+c}4}
=i^{-q}=(-1)^{\f{q+1}2}i$. For $r>2$ we have $r=2t\ge 4$, $2^4\mid
x+c$ and so $i^{-\f{x+c}4}=1$. If $q\e 1\mod 4$, from (4.3) and the
above we get
$$\align i^k&=\Ls 2q\Qs{1+i}q\Qs{A-Bi}q
\\&=(-1)^{\f{q-1}4}i^{\f{q-1}4} (-1)^{\f{x+1}4+\f
d4} i^{-\f{x+c}4+\f{q-1}4}\Qs x{c+di}\Qs{y^{-1}}{c+di}
\\&=(-1)^{\f{x+1}4+\f d4}i^{-\f{x+c}4}\Qs{x/y}{c+di}
=(-1)^{\f{x+c}4+\f d4+\f{1-c}4}i^{-\f{x+c}4}\Qs{x/y}{c+di}
\\&=(-1)^{\f{1-c}4}i^{-\f{x+c}4}\Qs{x/y}{c+di}
=(-1)^{\f{p-1}8}i^{-\f{x+c}4}\Qs{x/y}{c+di}.
\endalign$$
and so
$$\Qs{x/y}{c+di}= \cases (-1)^{\f{p-1}8}i^{k+1}
&\t{if $r=2$,}
\\(-1)^{\f{p-1}8}i^k&\t{if $r>2$.}\endcases$$
Now applying Lemma 2.8 we deduce the result.
 If $q\e 3\mod 4$, then
$\sls y{c^2+d^2}=\sls{y_0}p=\sls p{y_0}=\sls{x^2+2qy^2}{y_0}=1$.
From (4.3) and the above we deduce that
$$\align \Qs{A-Bi}q&=(-1)^{\f{x+1}4+\f d4}i^{-\f{x+c}4+\f{-q-1}4}
\Qs x{c-di}\Qs{y^{-1}}{c+di}
\\&=(-1)^{\f{x+c}4+\f d4+\f{1-c}4}i^{-\f{x+c}4-\f{q+1}4}
\Qs x{c+di}^{-1}\Qs y{c+di}\Qs{y^{-2}}{c+di}
\\&=(-1)^{\f{1-c}4}i^{-\f{x+c}4-\f{q+1}4}
\Qs{x^{-1}y}{c+di}\Qs{y^2}{c+di}^{-1}
\\&=(-1)^{\f{c^2-1}8}i^{-\f{x+c}4-\f{q+1}4}\Qs{x/y}{c+di}^{-1}\Ls
y{c^2+d^2}^{-1}
\\&=(-1)^{\f{p-1}8}i^{-\f{x+c}4-\f{q+1}4}\Qs{x/y}{c+di}^{-1}.
\endalign$$
Thus, using (4.3) we see that
$$i^k=\Qs{1+i}q\Qs{A-Bi}q^{-1}
=i^{\f{-q-1}4}\cdot
(-1)^{\f{p-1}8}i^{\f{x+c}4+\f{q+1}4}\Qs{x/y}{c+di}$$ and so
$$\Qs{x/y}{c+di}=(-1)^{\f{p-1}8}i^{k-\f{x+c}4}
=\cases (-1)^{\f{p-1}8}i^{k-q}=(-1)^{\f{p-1}8}i^{k+1}&\t{if $r=2$,}
\\(-1)^{\f{p-1}8}i^k&\t{if $r>2$.}
\endcases$$
Applying Lemma 2.8 we deduce the result.
 \par
 {\bf Case 3.} $m>r$. Using
Lemmas 2.1-2.5 we see that
$$\align i^k&=\Qs{-d/(x+c)+i}q=
\Qs{\f d{2^r}-\f{x+c}{2^r}i}q=(-1)^{\f{q-1}2\cdot\f{x+c}{2^{r+1}}}
\Qs q {\f d{2^r}-\f{x+c}{2^r}i}
\\&=(-1)^{\f{q-1}2\cdot\f{x+c}{2^{r+1}}}
\Qs {qy^2} {\f d{2^r}-\f{x+c}{2^r}i} \Qs {y^2} {\f
d{2^r}-\f{x+c}{2^r}i}
\\&=(-1)^{\f{q-1}2\cdot\f{x+c}{2^{r+1}}}
\Qs {((x+c)^2+d^2)/2-x(x+c)} {\f d{2^r}-\f{x+c}{2^r}i} \Ls y{\sls
d{2^r}^2+\sls{x+c}{2^r}^2}
\\&=(-1)^{\f{q-1}2\cdot \f{x+c}{2^{r+1}}}
\Qs {-2^mx(x+c)/2^m} {\f d{2^r}-\f{x+c}{2^r}i} \Ls
y{((x+c)^2+d^2)/2^{2r}}
\\&=(-1)^{\f{q-1}2\cdot \f{x+c}{2^{r+1}}}\Qs 2{\f
d{2^r}-\f{x+c}{2^r}i}^m\Qs {-x(x+c)/2^m} {\f d{2^r}-\f{x+c}{2^r}i}
\Ls{2^ty_0}{((x+c)^2+d^2)/2^{2r}}
\\&=(-1)^{\f{q-1}2\cdot \f{x+c}{2^{r+1}}}i^{-\f
d{2^r}\cdot\f{x+c}{2^{r+1}}\cdot
m}(-1)^{\f{-x(x+c)/2^m-1}2\cdot\f{x+c}{2^{r+1}}}\\&\q\times\Qs {\f
d{2^r}-\f{x+c}{2^r}i}{-x(x+c)/2^m}(-1)^{\f{x+c}{2^{r+1}}t}
\Ls{(x+c)^2+d^2}{y_0}
 .\endalign$$
 Applying  Lemma 2.7 we derive that
 $$\align &\Qs {\f
d{2^r}-\f{x+c}{2^r}i}{-x(x+c)/2^m} \Ls{(x+c)^2+d^2}{y_0}
\\&=\Qs{d-(x+c)i}{-x}\Qs{d-(x+c)i}{(x+c)/2^m}
\Ls{2qy^2+2x(c+x)}{y_0} \\&=\Qs{d-ci}{-x}\Qs
d{(x+c)/2^m}\Ls{2x(c+x)}{y_0}
\\&=\Qs{-i}x\Qs{c+di}x\Ls 2{y_0}\Qs{d+ci}{y_0}
=(-1)^{\f{x+1}4}\Qs x{c+di}\Qs {i(d+ci)}{y_0}
\\&=(-1)^{\f{x+1}4}\Qs x{c+di}\Qs{c-di}{y_0}
=(-1)^{\f{x+1}4}\Qs x{c+di}\Qs {y_0}{c-di}
\\&=(-1)^{\f{x+1}4}\Qs x{c+di}\Qs {y_0}{c+di}^{-1}
=(-1)^{\f{x+c}4+\f{1-c}4}\Qs x{c+di}\Qs {2^t/y}{c+di}
\\&=(-1)^{\f{c-1}4+\f d4t}\Qs{x/y}{c+di}=(-1)^{\f{p-1}8+\f
d4t}\Qs{x/y}{c+di}.\endalign$$ Therefore,
$$\aligned i^k&=(-1)^{\f{q-1}2\cdot \f{x+c}{2^{r+1}}}i^{-\f{x+c}{2^{r+1}}\cdot
m}(-1)^{\f{x(x+c)/2^m+1}2\cdot\f{x+c}{2^{r+1}}}\cdot
(-1)^{\f{x+c}{2^{r+1}}t+\f{p-1}8+\f d4t}\Qs{x/y}{c+di}
\\&=\cases (-1)^{\f
{q-1}2+\f{x(x+c)/2^{r+1}+1}2+\f{p-1}8+(1+\f
d4)t}i^{-\f{x+c}{2^{r+1}}(r+1)}\qs{x/y}{c+di} &\t{if $m=r+1$,}
\\(-1)^{m\f{x+c}{2^{r+2}}+\f{p-1}8+\f d4t}\qs{x/y}{c+di}&\t{if $m>r+1$.}
\endcases\endaligned$$

 Since $2^t\ \Vert\ y$, $2^r\ \Vert\ d$, $2^m\ \Vert\
x+c$, $m>r\ge 2$ and $2qy^2-d^2=(x+c)^2-2x(x+c)$, we see that
$2^{m+1}\ \Vert\ 2qy^2-d^2$. That is, $2^{m+1}\ \Vert\
2^{2t+1}qy_0^2-2^{2r}d_0^2.$ Hence $m+1=\t{min}\{2t+1,2r\}\le 2r$
and so $t>1$. If $m+1=2r$, then $t\ge r$ and
$$\aligned &x\f{x+c}{2^m}=2^{m-1}\Ls{x+c}{2^m}^2-2^{2t+1-2r}qy_0^2+d_0^2
\\&\e 1-2^{2t+1-2r}q\e\cases 1\mod 4&\t{if $t>r$,}
\\3\mod 4&\t{if $t=r$}.\endcases\endaligned$$
If $m+1<2r$, then $m=2t$, $t<r$ and
$$\aligned &x\f{x+c}{2^m}=2^{m-1}\Ls{x+c}{2^m}^2
-qy_0^2+2^{2r-1-2t}d_0^2
\\&\e -q+2^{2r-2t-1}\e\cases -q\mod 4&\t{if $t<r-1$,}
\\2-q\mod 4&\t{if $t=r-1$}.\endcases\endaligned$$
\par For $r=2$ we must have $m=3$ and so
$$\align i^k&=(-1)^{\f{q-1}2+\f{x(x+c)/8+1}2+\f{p-1}8+(1+\f d4)t}
i^{-\f{x+c}8\cdot 3}\Qs{x/y}{c+di}
\\&=(-1)^{\f{q-1}2+\f{x(x+c)/8+1}2+\f{p-1}8}
i^{\f{x+c}8}\Qs{x/y}{c+di} =(-1)^{\f{q-1}2+\f{p-1}8}i\Qs{x/y}{c+di}.
\endalign$$
Now assume $r>2$. If $m+1=2r$, then $m>r+1$ and so
$$i^k=(-1)^{\f{x+c}{2^{r+2}}+\f{p-1}8}\Qs{x/y}{c+di}
=(-1)^{2^{r-3}+\f{p-1}8}\Qs{x/y}{c+di}=(-1)^{\f
d8+\f{p-1}8}\Qs{x/y}{c+di}.$$ If $m+1<2r$ and $m=r+1$, then $t=\f
m2=\f{r+1}2$, $x\f{x+c}{2^{r+1}}\e -q+2^{2r-2t-1}=-q+2^{r-2}\e -q+\f
d4\mod 4$ and so
$$\align i^k&=(-1)^{\f{q-1}2+\f {x(x+c)/2^{r+1}+1}2+\f{p-1}8+t}
i^{-\f{x+c}{2^{r+1}}(r+1)}\Qs{x/y}{c+di}
\\&=(-1)^{\f{q-1}2+\f{1-q+\f d4}2+\f{p-1}8+\f{r+1}2}(-1)^{\f{r+1}2}
\Qs{x/y}{c+di} =(-1)^{\f{p-1}8+\f d8}\Qs{x/y}{c+di}.
\endalign$$
If $m+1<2r$ and $m>r+1$, then $r>3$, $m=2t$ and so
$i^k=(-1)^{\f{p-1}8}\qs{x/y}{c+di}=(-1)^{\f{p-1}8+\f
d8}\qs{x/y}{c+di}$. Hence
$$\Qs{x/y}{c+di}=\cases
(-1)^{\f{q-1}2+\f{p-1}8}i^{k-1}&\t{if $r=2$,}
\\(-1)^{\f{p-1}8+\f d8}i^k&\t{if $r>2$.}
\endcases$$
Now applying Lemma 2.8 we obtain the result.
\par By the above, the
theorem is proved.

  \pro{Theorem 4.3} Let
$p$ and $q$ be primes such that $p\e 1\mod 8$ and $q\e 3\mod 4$.
Suppose $p=c^2+d^2=x^2+2qy^2$, $c,d,x,y\in\Bbb Z$, $c\e 1\mod 4$,
$d=2^rd_0$, $d_0\e 1\mod 4$ and $\sls{c-di}x^{\f{q+1}4}\e i^m\mod
q$. Assume $(c,x+d)=1$ or $(d_0,x+c)=1$. Then
$$(-q)^{\f{p-1}8}\e (-1)^{\f {x-1}2\cdot \f{q+1}4}\Ls dc^m\mod p.$$
\endpro
Proof. Clearly $q\nmid x$ and $x$ is odd. We first assume
$(c,x+d)=1$. By the proof of Theorem 4.1,
$(q,(x+d)(c^2+(x+d)^2))=1$. It is easily seen that
$\f{c/(x+d)-i}{c/(x+d)+i}=\f{c-(x+d)i}{c+(x+d)i}\e \f{c-di}{ix}\mod
q$.
 Thus, for $k=0,1,2,3$, using Lemma 2.6 we get
$$\align &\Qs{c/(x+d)+i}q=i^k\\&\iff\f c{x+d}\in Q_k(q)
\iff \Ls{\f c{x+d}-i}{\f c{x+d}+i}^{\f{q+1}4}\e i^k\mod q
\\&\iff
\Ls{c-di}{ix}^{\f{q+1}4}\e i^k\mod q \iff \Ls{c-di}x^{\f{q+1}4}\e
i^{\f{q+1}4+k}\mod q.
\endalign$$
Since  $\sls{c-di}x^{\f{q+1}4}\e i^m\mod q$, from the above we
deduce that
$$\Qs{c/(x+d)+i}q=i^{m-\f{q+1}4}=\cases
(-1)^{\f{q+5}8}i^{m+1}&\t{if $q\e 3\mod
8$,}\\(-1)^{\f{q+1}8}i^m&\t{if $q\e 7\mod 8$.}\endcases$$ Now,
applying Theorem 4.1 we derive the result.
\par Now we assume $(d_0,x+c)=1$. By the proof of Theorem 4.2,
$(q,x+c)=(q,d^2+(x+c)^2)=1$.
 It is easily seen
that $\f{d+(x+c)i}{d-(x+c)i}\e \f{c-di}{-x}\mod q$.
 Thus, for $k=0,1,2,3$, using Lemma 2.6 we get
$$\align &\Qs{-d/(x+c)+i}q=i^k\\&\iff-\f d{x+c}\in Q_k(q)
\iff \Ls{-\f d{x+c}-i}{-\f d{x+c}+i}^{\f{q+1}4}\e i^k\mod q
\\&\iff \Ls{d+(x+c)i}{d-(x+c)i}^{\f{q+1}4}\e i^k\mod q\iff
\Ls{c-di}{-x}^{\f{q+1}4}\e i^k\mod q
\\&\iff \Ls{c-di}x^{\f{q+1}4}\e i^{\f{q+1}2+k}\mod q.
\endalign$$
Since  $\sls{c-di}x^{\f{q+1}4}\e i^m\mod q$, from the above we
deduce that $\sqs{-d/(x+c)+i}q=i^{m-\f{q+1}2}=(-1)^{\f{q+1}4}i^m$.
Now applying Theorem 4.2 we deduce the result. The proof is now
complete.

\pro{Corollary 4.1} Let $p$ and $q$ be primes such that $p\e 1\mod
8$ and $q\e 3\mod 8$. Suppose $p=c^2+d^2=x^2+2qy^2$, $c,d,x,y\in\Bbb
Z$, $c\e 1\mod 4$, $d=2^rd_0$, $q\mid cd$ and $d_0\e 1\mod 4$.
Assume $(c,x+d)=1$ or $(d_0,x+c)=1$. Then
$$(-q)^{\f{p-1}8}\e \cases \pm (-1)^{\f {x-1}2}\mod p&\t{if
$x\e \pm c\mod q$,}\\\mp(-1)^{\f{q-3}8+\f{x-1}2}\f dc\mod p&\t{if
$x\e \pm d\mod q$.}\endcases$$
\endpro
Proof. If $x\e \pm c\mod q$, then $q\mid d$ and so
$\sls{c-di}x^{\f{q+1}4}\e (\pm 1)^{\f{q+1}4}= \pm 1\mod q$. If $x\e
\pm d\mod q$, then $q\mid c$ and so $\sls{c-di}x^{\f{q+1}4}\e (\mp
i)^{\f{q+1}4}= \mp (-1)^{\f{q-3}8}i\mod q$. Now applying Theorem 4.3
we deduce the result.
\par As an example, taking $q=3$ in Corollary 4.1 we see that if
 $p$ is a prime of the form $24k+1$ and so
$p=c^2+d^2=x^2+6y^2$, then
$$(-3)^{\f{p-1}8}\e\cases \pm (-1)^{\f{x-1}2}\mod p&\t{if $x\e \pm
c\mod 3$,}
\\\mp (-1)^{\f{x-1}2}\f dc\mod p&\t{if $x\e \pm d\mod 3$.}
\endcases\tag 4.5$$

\pro{Theorem 4.4} Let $p$ and $q$ be primes such that $p\e 1\mod 8$,
$q\e 7\mod 8$, $p=c^2+d^2=x^2+2qy^2$, $c,d,x,y\in\Bbb Z$, $c\e 1\mod
4$, $d=2^rd_0$ and $d_0\e 1\mod 4$.  Assume $(c,x+d)=1$ or
$(d_0,x+c)=1$. Then
 $$(-q)^{\f{p-1}8}\e \Ls dc^m\mod p\iff
\Ls{c-di}{c+di}^{\f{q+1}8}\e i^m\mod q.$$
\endpro
Proof. Observe that
$$\Ls{c-di}{c+di}^{\f{q+1}8}=\f{(c-di)^{\f{q+1}4}}{(c^2+d^2)^{\f{q+1}8}}
=\f{(c-di)^{\f{q+1}4}}{(x^2+2qy^2)^{\f{q+1}8}}\e
\Ls{c-di}x^{\f{q+1}4}\mod q.$$ The result follows from Theorem 4.3.
\par\q
\par We note that if $q\nmid d$, then the $m$ in Theorem 4.4
depends only on $\f cd\mod q$.
 \pro{Corollary 4.2} Let $p\e 1\mod 8$ and $q\e 7\mod{8}$ be
primes such that $p=c^2+d^2=x^2+2qy^2$ with $c,d,x,y\in\Bbb Z$ and
$q\mid cd(c^2-d^2)$. Suppose $c\e 1\mod 4$, $d=2^rd_0$ and $d_0\e
1\mod 4$. Assume $(c,x+d)=1$ or $(d_0,x+c)=1$. Then
$$q^{\f{p-1}8}\e \cases (-1)^{\f{q+1}8}
\mod p&\t{if $q\mid c$,}
\\1\mod p&\t{if $q\mid d$,}
\\\pm (-1)^{\f{q+9}{16}}\f dc\mod p&\t{if $16\mid (q-7)$ and
$c\e \pm d\mod q$,}
\\(-1)^{\f{q+1}{16}}\mod p&\t{if $16\mid (q-15)$ and
$c\e \pm d\mod q$.}
\endcases$$
\endpro
Proof. Clearly
$$\f{c-di}{c+di}\e \cases -1\mod q&\t{if $q\mid c$,}
\\1\mod q&\t{if $q\mid d$,}
\\-i\mod q&\t{if $c\e d\mod q$,}
\\i\mod q&\t{if $c\e -d\mod q$.}
\endcases$$
Thus the result follows from Theorem 4.4. \pro{Theorem 4.5} Let $p$
and $q$ be distinct primes, $p\e 1\mod 8$, $q\e 1\mod 4$,
$p=c^2+d^2=x^2+2qy^2$, $q=a^2+b^2$, $a,b,c,d,x,y\in\Bbb Z$, $c\e
1\mod 4$, $d=2^rd_0$ and $d_0\e 1\mod 4$. Assume $(c,x+d)=1$ or
$(d_0,x+c)=1$. Suppose $\sls{ac+bd}{ax}^{\f{q-1}4}\e \sls ba^m\mod
q$. Then
$$(-q)^{\f{p-1}8}
\e(-1)^{\f {x-1}2\cdot\f{q-1}4+\f d4+\f y2}\Ls dc^m\mod p.$$
\endpro
Proof. Clearly $q\nmid x$. We first assume $(c,x+d)=1$. By the proof
of Theorem 4.1, $(q,(x+d)(c^2+(x+d)^2))=1$. It is easily seen that
$\f{ac+b(x+d)}{ac-b(x+d)}\e \f{ac+bd}{ax}\cdot\f ba\mod q$.
 Thus, for $k=0,1,2,3$, using Lemma 2.6 we get
$$\align \Qs{c/(x+d)+i}q=i^k&\iff\f c{x+d}\in Q_k(q)
\iff \Ls{\f c{x+d}+\f ba}{\f c{x+d}-\f ba}^{\f{q-1}4}\e \Big(\f
ba\Big)^k\mod q
\\&\iff \Ls{ac+b(x+d)}{ac-b(x+d)}^{\f{q-1}4}\e \Ls ba^k\mod q
\\&\iff
\Big(\f{ac+bd}{ax}\cdot\f ba\Big)^{\f{q-1}4}\e \Ls ba^k\mod q
\\&\iff \Ls{ac+bd}{ax}^{\f{q-1}4}\e \Ls ba^{k-\f{q-1}4}\mod q.
\endalign$$
Since $ \sls{ac+bd}{ax}^{\f{q-1}4}\e \sls ba^m\mod q$, from the
above we get $\sqs{c/(x+d)+i}q=i^{m+\f{q-1}4}$.  Now the result
follows from Theorem 4.1 immediately.
\par Suppose $(d_0,x+c)=1$. By the proof of Theorem 4.2, $(q,(x+c)(d^2+(x+c)^2))=1$. It is easily seen that
$\f{ad-b(x+c)}{ad+b(x+c)}\e \f{ac+bd}{-ax}\mod q$.
 Thus, for $k=0,1,2,3$, using Lemma 2.6 we get
$$\align &\Qs{-d/(x+c)+i}q=i^k\\&\iff-\f d{x+c}\in Q_k(q)
\iff \Ls{-\f d{x+c}+\f ba}{-\f d{x+c}-\f ba}^{\f{q-1}4}\e \Big(\f
ba\Big)^k\mod q
\\&\iff \Ls{ad-b(x+c)}{ad+b(x+c)}^{\f{q-1}4}\e \Ls ba^k\mod q
\iff \Ls{ac+bd}{-ax}^{\f{q-1}4}\e \Ls ba^k\mod q
\\&\iff \Ls{ac+bd}{ax}^{\f{q-1}4}\e \Ls ba^{\f{q-1}2+k}\mod q.
\endalign$$
Since $ \sls{ac+bd}{ax}^{\f{q-1}4}\e \sls ba^m\mod q$, by the above
we get $\sqs{-d/(x+c)+i}q=i^{m-\f{q-1}2}$. Now applying Theorem 4.2
  we derive the result. The proof is now
complete.

\pro{Corollary 4.3} Let $p\e 1\mod 8$ and $q\e 5\mod 8$ be primes
such that $p=c^2+d^2=x^2+2qy^2$ with $c,d,x,y\in\Bbb Z$ and $q\mid
cd$. Suppose $c\e 1\mod 4$, $d=2^rd_0$ and $d_0\e 1\mod 4$. Assume
$(c,x+d)=1$ or $(d_0,x+c)=1$. Then
$$(-q)^{\f{p-1}8}\e \cases
\pm(-1)^{\f d4+\f{x-1}2+\f y2}\mod p&\t{if $x\e\pm c\mod{q}$,}
\\\pm (-1)^{\f{q-5}8+\f d4+\f{x-1}2+\f y2}\f dc\mod p
&\t{if $x\e \pm d\mod{q}$.}
\endcases$$
\endpro
Proof. Suppose $q=a^2+b^2$ with $a,b\in\Bbb Z$. If $x\e \pm c\mod
q$, then $q\mid d$ and so $\sls{ac+bd}{ax}^{\f{q-1}4}\e \sls
cx^{\f{q-1}4}\e (\pm 1)^{\f{q-1}4} = \pm 1\mod q$. If $x\e \pm d\mod
q$, then $q\mid c$ and so $\sls{ac+bd}{ax}^{\f{q-1}4}\e \sls
{bd}{ax}^{\f{q-1}4}\e (\pm \f ba)^{\f{q-1}4} \e \pm
(-1)^{\f{q-5}8}\f ba\mod q$. Now putting the above with Theorem 4.5
we deduce the result.

\pro{Theorem 4.6} Let $p$ and $q$ be distinct primes such that $p\e
1\mod 8$, $q\e 1\mod 8$, $p=c^2+d^2=x^2+2qy^2$, $q=a^2+b^2$,
$a,b,c,d,x,y\in\Bbb Z$, $c\e 1\mod 4$, $d=2^rd_0$ and $d_0\e 1\mod
4$. Assume $(c,x+d)=1$ or $(d_0,x+c)=1$.  Then
 $$(-q)^{\f{p-1}8}\e (-1)^{\f d4+\f y2}\Ls dc^m\mod p\iff
 \Ls{ac+bd}{ac-bd}^{\f{q-1}8}\e \Ls ba^m\mod q.$$
\endpro

Proof. Observe that $b^2\e -a^2\mod q$, $p\e x^2\mod q$ and so
$$\Ls{ac+bd}{ac-bd}^{\f{q-1}8}=\f{(ac+bd)^{\f{q-1}4}}{(a^2c^2-b^2d^2
)^{\f{q-1}8}} \e\f{(ac+bd)^{\f{q-1}4}}{(a^2p)^{\f{q-1}8}}\e
\Ls{ac+bd}{ax}^{\f{q-1}4}\mod q.$$ The result follows from Theorem
4.5.
\par\q
\par We note that if $q\nmid d$, then the $m$ in Theorem 4.6
depends only on $\f cd\mod q$. \pro{Corollary 4.4} Let $p$ and $q$
be distinct primes of the form $8k+1$ such that
$p=c^2+d^2=x^2+2qy^2$ with $c,d,x,y\in\Bbb Z$ and $q\mid
cd(c^2-d^2)$. Suppose $c\e 1\mod 4$, $d=2^rd_0$ and $d_0\e 1\mod 4$.
Assume $(c,x+d)=1$ or $(d_0,x+c)=1$. Then
$$(-q)^{\f{p-1}8}\e \cases (-1)^{\f{q-1}8+\f d4+\f y2}
\mod p&\t{if $q\mid c$,}
\\(-1)^{\f d4+\f y2}\mod p&\t{if $q\mid d$,}
\\ (-1)^{\f{q-1}{16}+\f d4+\f y2}\mod p
&\t{if $16\mid (q-1)$ and $c\e \pm d\mod q$,}
\\\pm (-1)^{\f{q-9}{16}+\f d4+\f y2}\f dc\mod p
&\t{if $16\mid (q-9)$ and $c\e \pm d\mod q$.}
\endcases$$
\endpro
Proof. Suppose that $q=a^2+b^2$ with $a,b\in\Bbb Z$. Then clearly
$$\Ls{ac+bd}{ac-bd}^{\f{q-1}8}\e \cases (-1)^{\f{q-1}8}\mod q&\t{if $q\mid c$,}
\\1\mod q&\t{if $q\mid d$,}
\\(-1)^{\f{q-1}{16}}\mod q&\t{if $16\mid (q-1)$ and $c\e \pm d\mod q$,}
\\\pm (-1)^{\f{q-9}{16}}\f ba\mod q&\t{if $16\mid (q-9)$ and $c\e \pm d\mod q$.}
\endcases$$
Thus the result follows from Theorem 4.6.

\pro{Theorem 4.7} Let $p\e 1\mod 8$ be a prime,
$p=c^2+d^2=x^2+2(a^2+b^2)y^2$, $a,b,c,d,x,y\in\Bbb Z$, $a\not=0$,
$2\mid a$, $(a,b)=1$, $c\e 1\mod 4$, $d=2^rd_0$ and $d_0\e 1\mod 4$.
Assume $(c,x+d)=1$ or $(d_0,x+c)=1$. Then
$$\aligned &(-a^2-b^2)^{\f{p-1}8}\e\cases
(-1)^{\f d4+\f y2}\sls cd^m\mod p&\t{if $4\mid a$,}
\\(-1)^{\f{b-1}2+\f d4+\f y2+\f{x-1}2}\sls cd^{m-1}\mod p
&\t{if $4\mid a-2$.}
\endcases\\&\iff \Qs{(ac+bd)/x}{b+ai}=i^m.\endaligned$$
\endpro
Proof. Suppose $q=a^2+b^2$ and $\sqs{(ac+bd)/x}{b+ai}=i^m$.  Then
clearly $q\e 1\mod 4$ and $p\nmid q$.
 We first assume $(c,x+d)=1$.
By the proof of Theorem 4.1, $(q,x+d)=(q,c^2+(x+d)^2)=1$. Since
$\f{c-(x+d)i}{c+(x+d)i}\e \f{c-di}{ix}\mod q$,  from [S5, p.24] we
know that
$$\Qs{c/(x+d)+i}q=(-1)^{\f {b+1}2\cdot \f a2+[\f
q8]}i^{-m}.$$ This together with Theorem 4.1 yields the result in
this case.
\par Now we assume $(d_0,x+c)=1$. By the proof of Theorem 4.2,
$(q,x+c)=(q,(x+c)^2+d^2)=1$. Since $\f{d+(x+c)i}{d-(x+c)i}\e
\f{c-di}{-x}\mod q$, from [S5, p.24] we know that
$$\Qs{-d/(x+c)+i}q=\cases (-1)^{\f{b+1}2}i^{1-m}
&\t{if $4\mid a-2$,}\\i^{-m}&\t{if $4\mid a$.}
\endcases$$
Now applying Theorem 4.2 we deduce the result in this case. So the
theorem is proved. \pro{Corollary 4.5} Let $p\e 1,9\mod{40}$ be a
prime and so $p=c^2+d^2=x^2+10y^2$ with $c,d,x,y\in\Bbb Z$. Suppose
$c\e 1\mod 4$, $d=2^rd_0$ and $d_0\e 1\mod 4$. Assume $(c,x+d)=1$ or
$(d_0,x+c)=1$. Then
$$(-5)^{\f{p-1}8}\e\cases \pm (-1)^{\f d4+\f {x-1}2+\f y2}\mod
p&\t{if $\f{2c-d}x\e \pm 2\mod 5$,}
\\\pm (-1)^{\f d4+\f{x-1}2+\f y2}\f cd\mod p&\t{if $\f{2c-d}x\e \pm
1\mod 5$.}
\endcases$$
\endpro
Proof. Since $\sqs{\pm 1}{1+2i}=\pm 1$ and $\sqs{\pm 2}{1+2i}=\pm
i$, taking $a=2$ and $b=1$ in Theorem 4.7 we derive the result.
\par We remark that Corollary 4.5 partially solves [S4, Conjecture
9.8].
 \subheading{5. Congruences for $U_{\f{p-1}4}(2a,-1)$ and
$V_{\f{p-1}4}(2a,-1)\mod p$}
\par For two numbers $P$ and $Q$ the Lucas sequences $\{U_n(P,Q)\}$ and $\{V_n(P,Q)\}$
are defined by
$$\align&U_0(P,Q)=0,\ U_1(P,Q)=1,\ U_{n+1}(P,Q)=PU_n(P,Q)
-QU_{n-1}(P,Q)\ (n\ge 1), \\&V_0(P,Q)=2,\ V_1(P,Q)=P,\
V_{n+1}(P,Q)=PV_n(P,Q)- QV_{n-1}(P,Q)\ (n\ge 1).\endalign$$ Set
$D=P^2-4Q$. It is well known that
$$U_n(P,Q)=\f 1{\sqrt{D}}\Big\{\Big(\f{P+\sqrt D}2\Big)^n-
\Big(\f{P-\sqrt{D}}2\Big)^n\Big\}\q(D\not=0)\tag 5.1$$ and
$$V_n(P,Q)=\Big(\f{P+\sqrt{D}}2\Big)^n+
\Big(\f{P-\sqrt{D}}2\Big)^n.\tag 5.2$$
 \pro{Theorem 5.1} Let $p$ be a prime of
the form $8k+1$ and $a\in\Bbb Z$ with $2\nmid a$. Suppose that
$p=c^2+d^2=x^2+(a^2+1)y^2$, $c,d,x,y\in\Bbb Z$, $c\e 1\mod 4$,
$d=2^rd_0$, $y=2^{\beta}y_0$ and $d_0\e y_0\e 1\mod 4$. Assume
$(c,x+d)=1$ or $(d_0,x+c)=1$. Then
$$U_{\f{p-1}4}(2a,-1)\e\cases (-1)^{\f{a-1}2+\f d4+\f{x-1}2}\f yx\mod p&\t{if
$4\mid y-2$,}\\0\mod p&\t{if $4\mid y$}
\endcases$$ and
$$V_{\f{p-1}4}(2a,-1)\e\cases 0\mod p&\t{if
$4\mid y-2$,}\\2(-1)^{\f d4+\f y4}\mod p&\t{if $4\mid y$.}
\endcases$$
\endpro
Proof. Set $a_1=(1-(-1)^{\f{a-1}2}a)/2$ and
$b_1=(1+(-1)^{\f{a-1}2}a)/2$. Then $2\mid a_1$, $2\nmid b_1$ and
$a^2+1=2(a_1^2+b_1^2)$.  It is clear that
$$\Qs{(a_1c+b_1d)/((-1)^{\f{x-1}2}x)}{b_1+a_1i}
=(-1)^{\f{x-1}2\cdot\f{a_1}2}\Qs{(a_1c+b_1d)/x} {b_1+a_1i}.$$
 We first assume $a\e 1\mod 4$. Replacing $d,x$ with $-d,
 (-1)^{\f{x-1}2}x$ in [S4, Theorem 8.3(i)] we obtain
$$U_{\f{p-1}4}(2a,-1)\e\cases \mp(-1)^{\f{x-1}2\cdot\f{a-1}4}(-a_1^2-b_1^2)^{\f{p-1}8}(-\f
cd)^{(1-(-1)^{\f{a-1}4})/2}(-1)^{\f{x-1}2}\f yx\mod p\\\qq\qq\q\t{if
$4\mid y-2$ and $\qs{(a_1c+b_1d)/x}{b_1+a_1i}=\pm 1$,}
\\\mp (-1)^{\f{x-1}2\cdot\f{a-1}4}(-a_1^2-b_1^2)^{\f{p-1}8}(-\f
cd)^{1+(1-(-1)^{\f{a-1}4})/2}(-1)^{\f{x-1}2}\f yx\mod p
\\\qq\qq\q\t{if $4\mid y-2$ and
$\qs{(a_1c+b_1d)/x}{b_1+a_1i}=\pm i$,}
\\0\mod p\ \;\t{if $4\mid y$}
\endcases$$
and
$$V_{\f{p-1}4}(2a,-1)\e\cases \pm 2(-1)^{\f{x-1}2\cdot\f{a-1}4+\f y4}(-a_1^2-b_1^2)^{\f{p-1}8}(-\f
cd)^{(1-(-1)^{\f{a-1}4})/2}\mod p\\\qq\qq\q\t{if $4\mid y$ and
$\qs{(a_1c+b_1d)/x}{b_1+a_1i}=\pm 1$,}
\\\pm2(-1)^{\f{x-1}2\cdot\f{a-1}4+\f y4}(-a_1^2-b_1^2)^{\f{p-1}8}(-\f
cd)^{1+(1-(-1)^{\f{a-1}4})/2}\mod p
\\\qq\qq\q\t{if $4\mid y$ and
$\qs{(a_1c+b_1d)/x}{b_1+a_1i}=\pm i$,}
\\0\mod p\ \;\t{if $4\mid y-2$.}
\endcases$$
 From Theorem 4.7 we know that
$$\aligned &(-a_1^2-b_1^2)^{\f{p-1}8}
\\&\e\cases \pm (-1)^{\f d4+\f y2}\mod p&\t{if $4\mid a_1$ and
$\sqs{(a_1c+b_1d)/x}{b_1+a_1i}=\pm 1$,}
\\\pm (-1)^{\f{b_1-1}2+\f d4+\f y2+\f{x-1}2}\f dc\mod p&\t{if $4\mid
a_1-2$ and $\sqs{(a_1c+b_1d)/x}{b_1+a_1i}=\pm 1$,}
\\\pm (-1)^{\f d4+\f y2}\f cd\mod p&\t{if $4\mid a_1$ and
$\sqs{(a_1c+b_1d)/x}{b_1+a_1i}=\pm i$,}
\\\pm(-1)^{\f{b_1-1}2+\f d4+\f y2+\f{x-1}2}\mod p
&\t{if $4\mid a_1-2$ and $\sqs{(a_1c+b_1d)/x}{b_1+a_1i}=\pm i$.}
\endcases\endaligned$$
Now putting the above together we deduce the result in the case $a\e
1\mod 4$. The case $a\e 3\mod 4$ can be proved similarly by using
[S4, Theorem 8.3(ii)] and Theorem 4.7. \pro{Corollary 5.1} Let $p$
be a prime of the form $8k+1$ and $a\in\Bbb Z$ with $2\nmid a$.
Suppose that $p=c^2+d^2=x^2+(a^2+1)y^2$, $c,d,x,y\in\Bbb Z$, $c\e
1\mod 4$, $d=2^rd_0$, $d_0\e  1\mod 4$ and $4\mid y$. Assume
$(c,x+d)=1$ or $(d_0,x+c)=1$. Then
$$(a+\sqrt{a^2+1})^{\f{p-1}4}\e (-1)^{\f d4+\f y4}\mod p.$$
\endpro
Proof. From (5.1), (5.2) and Theorem 5.1 we see that
$$\align (a+\sqrt{a^2+1})^{\f{p-1}4}&=\f
12V_{\f{p-1}4}(2a,-1)+\sqrt{a^2+1}U_{\f{p-1}4}(2a,-1)
\\&\e \f
12V_{\f{p-1}4}(2a,-1)\e (-1)^{\f d4+\f y4}\mod p.\endalign$$ This
proves the corollary.
\pro{Corollary 5.2} Let $p$ be a prime of the
form $8k+1$ and $a\in\Bbb Z$ with $2\nmid a$. Suppose that
$p=c^2+d^2=x^2+(a^2+1)y^2$, $c,d,x,y\in\Bbb Z$, $c\e 1\mod 4$,
$d=2^rd_0$ and $d_0\e 1\mod 4$. Assume $(c,x+d)=1$ or $(d_0,x+c)=1$.
Then
$$p\mid U_{\f{p-1}8}(2a,-1)\iff 4\mid y\qtq{and}
\f{p-1}8\e\f d4+\f y4\mod 2.$$
\endpro
Proof. By [S4, (1.5)],
$$p\mid U_{\f{p-1}8}(2a,-1)\iff V_{\f{p-1}4}(2a,-1)\e
2(-1)^{\f{p-1}8}\mod p.$$ Now applying Theorem 5.1 we obtain the
result.\par\q
\newline{\bf Remark 5.1} Theorem 5.1 and Corollary 5.2 were
conjectured by the author in [S4, Conjectures 9.17 and 9.19].

\enddocument
\bye